\documentclass[fullpage, oneside]{amsart}

\usepackage{latexsym}
\usepackage{accents}
\usepackage[english]{babel}
\usepackage{hyperref}
\usepackage{a4wide}
\usepackage{amsthm,amssymb,amsmath,a4wide,graphicx,tikz,mathtools}
\usepackage{subfigure}
\usepackage{caption}
\usepackage{hhline}
\usetikzlibrary{patterns}
\usepackage{float}
\usepackage{dsfont}
\usepackage{mathtools}

\usepackage{titlesec}
\usepackage{pgf}
\titlelabel{\thetitle.\quad}

\usepackage{siunitx}
\usetikzlibrary{arrows,shapes,automata,positioning}

\usepackage[a4paper]{geometry}
\makeatletter
\renewcommand*\env@matrix[1][*\c@MaxMatrixCols c]{%
  \hskip -\arraycolsep
  \let\@ifnextchar\new@ifnextchar
  \array{#1}}
\makeatother

\definecolor{apricot}{rgb}{0.98, 0.81, 0.69}
\definecolor{aquamarine}{rgb}{0.5, 1.0, 0.83}
\definecolor{babyblueeyes}{rgb}{0.63, 0.79, 0.95}
\definecolor{bananamania}{rgb}{0.98, 0.91, 0.71}
\definecolor{bittersweet}{rgb}{1.0, 0.44, 0.37}
\definecolor{bluebell}{rgb}{0.64, 0.64, 0.82}
\definecolor{blush}{rgb}{0.87, 0.36, 0.51}

\definecolor{cerulean}{rgb}{0.0, 0.48, 0.65}
\definecolor{darkcyan}{rgb}{0.0, 0.55, 0.55}

\definecolor{airforceblue}{rgb}{0.36, 0.54, 0.66}
\definecolor{antiquefuchsia}{rgb}{0.57, 0.36, 0.51}
\definecolor{asparagus}{rgb}{0.53, 0.66, 0.42}
\definecolor{ballblue}{rgb}{0.13, 0.67, 0.8}
\definecolor{blueviolet}{rgb}{0.54, 0.17, 0.89}
\definecolor{brightgreen}{rgb}{0.4, 1.0, 0.0}
\definecolor{brightcerulean}{rgb}{0.11, 0.67, 0.84}
\definecolor{brightpink}{rgb}{1.0, 0.0, 0.5}
\definecolor{brightturquoise}{rgb}{0.03, 0.91, 0.87}
\definecolor{brightube}{rgb}{0.82, 0.62, 0.91}
\definecolor{brilliantlavender}{rgb}{0.98, 0.86, 1.0}
\definecolor{bluebell}{rgb}{0.64, 0.64, 0.82}
\definecolor{bluegray}{rgb}{0.4, 0.6, 0.8}
\definecolor{bluegreen}{rgb}{0.0, 0.87, 0.87}
\definecolor{brandeisblue}{rgb}{0.0, 0.44, 1.0}
\definecolor{capri}{rgb}{0.0, 0.75, 1.0}
\definecolor{fandango}{rgb}{0.71, 0.2, 0.54}
\definecolor{ashgrey}{rgb}{0.7, 0.75, 0.71}
\definecolor{gainsboro}{rgb}{0.86, 0.86, 0.86}
\definecolor{bubblegum}{rgb}{0.99, 0.76, 0.8}
\definecolor{brilliantrose}{rgb}{1.0, 0.33, 0.64}
\definecolor{brightmaroon}{rgb}{0.76, 0.13, 0.28}

\DeclareMathOperator{\KI}{KI}

\newtheorem{thm}{Theorem}[section]
\newtheorem{lem}[thm]{Lemma}

\newtheorem{defn}[thm]{Definition}
\newtheorem{ex}[thm]{Example}
\newtheorem{prop}[thm]{Proposition}
\newtheorem{nrem}[thm]{Remark}

\begin{document}
\title{Matching for random systems with an application to minimal weight expansions}
\author[Karma Dajani] {Karma Dajani$^\dagger$}
\author[Charlene Kalle]{Charlene Kalle$^\ddagger$}
\author[Marta Maggioni]{Marta Maggioni$^\ddagger$}

\address[$\dagger$]{Department of Mathematics, Utrecht University, P.O. Box 80010, 3508TA Utrecht, The Netherlands}
\email[Karma Dajani]{k.dajani1@uu.nl}
\address[$\ddagger$]{Mathematisch Instituut, Leiden University, Niels Bohrweg 1, 2333CA Leiden, The Netherlands}
\email[Charlene Kalle]{kallecccj@math.leidenuniv.nl}
\email[Marta Maggioni]{m.maggioni@math.leidenuniv.nl}


\begin{abstract}
We extend the notion of matching for one-dimensional dynamical systems to random matching for random dynamical systems on an interval. We prove that for a large family of piecewise affine random systems of the interval the property of random matching implies that any invariant density is piecewise constant. We further introduce a one-parameter family of random dynamical systems that produce signed binary expansions of numbers in the interval $[-1,1]$. This family has random matching for Lebesgue almost every parameter. We use this to prove that the frequency of the digit 0 in the associated signed binary expansions never exceeds $\frac12$.
\end{abstract}
\maketitle

\section{Introduction}
Optimal algorithms for the computation of powers of elements in a group are at the basis of many public key cryptosystems. Here the group is either the multiplicative group of a finite field or the group of points on an elliptic curve and the optimality refers to the ability of computing high powers in a short amount of time. One such algorithm is the {\em binary method}, introduced in \cite{Kn} and based on the binary expansion of the power. More precisely, if $x$ is an element of a given group, and $a = \sum_{k=0}^n d_k 2^k \in \mathbb{N}$ for some digits $d_k \in \{0,1\}$, then
\[x^a= \prod_{k=0}^n x^{d_k 2^k},\]
and the power $x^a$ is computed by taking the product of repeated squarings. While the number of squarings is given by the length $n$ of the binary expansion of $a$, the number of multiplications equals the number of non-zero bits $d_k$ in the expansion or its \textit{Hamming weight}. Clearly, a lower Hamming weight implies fewer multiplications and a faster result. To increase the number of zero bits, \cite{Boo} introduced a {\em signed binary representation}, i.e., a binary representation with digits in the set $\{-1, 0, 1\}$. This signed binary representation was later adopted in several methods in elliptic cryptosystems, see e.g.~\cite{CMO, HP06} and the references therein.
\vskip .2cm

The ordinary binary representation of an integer $a$ is uniquely determined, but this is not the case for the signed one. In fact, each integer has infinitely many signed binary representations, which led to the study of algorithms that choose the ones with minimal Hamming weight (see e.g.~\cite{MO90,KT,LK}). 
Typically a number has several signed binary representations with minimal weight (see \cite{GH}), but already in the 1960's Reitwiesner proved in \cite{Rei} that the signed representation is unique when adding the constraint $d_k d_{k+1}=0$. Such representations are known as \textit{signed separated binary expansions}, or SSB for short. In \cite{DKL06} it is shown how to obtain SSB expansions through the binary odometer and a three state Markov chain. Furthermore, in \cite{DKL06} the set $K:=\{d_1d_2 \ldots \in \{-1,0,1\}^{\mathbb{N}}\, : \,  \forall k \in \mathbb{N}, \, d_k d_{k+1}=0\}$ is introduced as a compactification of $\mathbb{Z}$. The authors identified $K$, endowed with the left shift $\sigma$, with the map $S(x)= 2x \mod \mathbb{Z}$ on the interval $[-\frac23, \frac23]$ through the conjugation
\begin{equation}\nonumber
\psi(d_1d_2 \ldots)= \sum_{k=1}^{\infty} \frac{d_k}{2^{k}}.
\end{equation}
This dynamical viewpoint allowed them to obtain metric properties of the system $(K, \sigma)$, such as a $\sigma$-invariant measure, the maximal entropy and the frequency of 0 in typical expansions. 
\vskip .2cm

In \cite{DK} this dynamical approach was further developed by considering a family of \textit{symmetric doubling maps} $\{S_{\alpha}: [-1,1] \rightarrow [-1,1] \}_{\alpha \in [1,2]}$ defined by $S_\alpha(x)= 2x - d \alpha$ and
\[d = \begin{cases}
-1, & \text{if } x \in [-1, -\frac12), \\
0, & \text{if } x \in [-\frac12, \frac12], \\
1, & \text{if } x \in (\frac12,1].
\end{cases}\]
The map $S$ from \cite{DKL06}, producing SSB expansions, is then easily identified with the map $S_{\frac32}$. For each $\alpha \in [1,2]$ iterations of $S_{\alpha}$ give a signed binary expansion of the form $x = \sum_{k=1}^{\infty} \frac{d_k}{2^{k}}$ with $d_k \in \{-1,0,1\}$ for each number $x \in [-1,1]$. The authors of \cite{DK} showed that the frequency of $0$ in such expansions depends continuously on the parameter $\alpha$ and takes its maximal value $\frac23$, corresponding to the minimal Hamming weight of $\frac13$, precisely for $\alpha \in [\frac65, \frac32]$. It follows that typically only $\frac13$ of the digits in the SSB expansions of integers is different from 0. The results from \cite{DK} are obtained by finding a detailed description of the unique invariant probability density $f_{\alpha}$ of $S_\alpha$ for each value $\alpha$ and then explicitly computing the frequency of the digit 0 using Birkhoff's Ergodic Theorem. The fact that the family $\{S_{\alpha}\}$ exhibits the dynamical phenomenon of matching was essential for these results.

\vskip .2cm
In this article we consider signed binary expansions in the framework of random dynamical systems. The advantage of random systems in this context is that a single random system produces many more number expansions per number than a deterministic map, allowing one to study the properties of many expansions simultaneously. See e.g.~\cite{DaKr,DaVr1,DaVr2,DK07,KaKe,DO} for the use of random systems in the study of different types of number expansions. We will introduce a family of random systems $\{R_{\alpha}\}_{\alpha \in [1,2]}$, called {\em random symmetric doubling maps}, such that each element $R_{\alpha}$ produces for typical numbers in the interval $[-1,1]$ infinitely many different signed binary expansions. This is contrary to the map $S_\alpha$, which produces a unique signed binary expansion for each number in $[-1,1]$. Our main result for the family $\{R_{\alpha}\}_{\alpha \in [1,2]}$ is that the frequency of the digit 0 in typical signed binary expansions produced by any of the maps $R_\alpha$ is at most $\frac12$, and therefore the Hamming weight is at least $\frac12$. This reinforces the result from \cite{DK} that the maps $S_\alpha$ with $\alpha \in \big[ \frac65, \frac32 \big]$ perform best in terms of minimal weight.

\vskip .2cm
We obtain this result from Birkhoff's Ergodic Theorem after gathering detailed knowledge on the invariant probability densities of the random maps $R_\alpha$. We first express these densities as infinite sums of indicator functions using the algebraic procedure from \cite{KM}. To compute the frequency of 0 we need to evaluate the Lebesgue integral of these densities over part of the domain and therefore we convert the infinite sums into finite sums. For this we introduce a random version of the dynamical concept of matching that is available for one-dimensional systems (see e.g.~\cite{NN08,DKS09,BCIT13,BORG13,BCK,BCMP,CIT18,CM18,KLMM}). Our definition of random matching properly extends the one-dimensional notion of matching and we illustrate the concept with examples of random continued fraction maps and random generalised $\beta$-transformations. We show that under mild certain conditions, if a random system of piecewise affine maps defined on the same interval has random matching, then any invariant probability density of the system is piecewise constant. The precise formulation of this statement and the conditions are given in the next section. Finally, we use this random matching property to show that for Lebesgue almost all parameters $\alpha$ the invariant density of the random systems $R_\alpha$, producing signed binary expansions, is in fact piecewise constant.

\vskip .2cm

The article is outlined as follows. The second section is devoted to random matching for random systems defined on an interval. We first recall some preliminaries on invariant measures for random interval maps. We then define the notion of random matching and state and prove the result about densities of random systems of piecewise affine maps with matching. We also discuss the examples of random continued fraction transformations and random generalised $\beta$-transformations. In the third section we introduce and discuss the family  $\{R_{\alpha}\}$ of random symmetric doubling maps and the corresponding signed binary expansions. We prove that $R_{\alpha}$ has random matching for Lebesgue almost all $\alpha \in [1,2]$. We also provide a full description of the {\em matching intervals}, i.e., intervals of parameters that exhibit comparable matching behaviour, and describe the invariant densities of the maps $R_\alpha$. Finally we prove that typically the frequency of the digit $0$ in the signed binary expansions produced by $R_\alpha$ does not exceed $\frac12$ for any parameter $\alpha$.

\section{Random matching}

\subsection{The definition of random matching}
Matching is a dynamical phenomenon observed in certain families of piecewise smooth interval maps. If $T: I \to I$ is such a map (so the domain $I$ is an interval of real numbers), then we say that $T$ has {\em matching} if for every discontinuity point $c$ of $T$ or of the derivative $T'$ the orbits of the left and right limits $T^k(c^-) = \lim_{x \uparrow c} T^k(x)$ and $T^k(c^+) = \lim_{x \downarrow c} T^k(x)$ eventually meet, i.e., if for each $c$ there exist positive constants $M=M_c$ and $Q=Q_c$, such that
\begin{equation}\label{q:onedmatching1}
T^M(c^-) = T^Q (c^+).
\end{equation}
$T$ is then said to have {\em strong matching} if, moreover, the orbits of the left and right limits have equal one-sided derivatives at the moment they meet, i.e., if besides \eqref{q:onedmatching1} it also holds that
\begin{equation}\label{q:onedmatching2}
(T^M)'(c^-) = (T^Q)'(c^+).
\end{equation}

\vskip .2cm
It was proven in \cite[Theorem 1.2]{BCMP} (see also Remark 1.3 in \cite{BCMP}) that for any piecewise smooth $T$ with strong matching, any invariant probability measure $\mu$ that is absolutely continuous with respect to the Lebesgue measure has a piecewise smooth density. For continued fraction transformations (as in \cite{NN08,DKS09,KLMM} for example) it seems that matching is sufficient to guarantee the existence of a piecewise smooth density (since this is sufficient to construct a natural extension with finitely many pieces). The strong matching condition then enforces some stability in the matching behaviour of certain one-parameter families of continued fraction maps, which becomes visible in the appearance of so called {\em matching intervals} in the parameter space: If such a family has strong matching for one parameter, then one can find an interval of parameters around it, such that all the corresponding transformations have matching in the same number of steps and with comparable orbits.

\vskip .2cm
In this section we extend the above definitions of matching and strong matching to random dynamical systems. With a random map we mean a system evolving in discrete time units in which at each step one of a number of transformations is chosen at random and applied. One way to describe a random map is with a pseudo-skew product transformation as follows. Let $\Omega \subseteq \mathbb N$ be the index set of the available maps, so we have a collection of transformations $\{ T_j : I \rightarrow I \}_{j \in \Omega}$ defined on the same interval $I$ at our disposal. Let $\sigma: \Omega^{\mathbb N} \to \Omega^{\mathbb N}$ be the left shift on one-sided sequences. The {\em random map} or {\em pseudo-skew product} $R: \Omega^{\mathbb N} \times I \to \Omega^{\mathbb N} \times I$ is defined by 
\[R(\omega,x) = (\sigma(\omega), T_{\omega_1}x).\]
So, the coordinates of $\omega$ determine which of the maps $T_j$ is applied at each time step. Let $\mathbf p = (p_j)_{j \in \Omega}$ be a positive probability vector, i.e., $p_j > 0$ for all $j \in \Omega$ and $\sum_{j \in \Omega} p_j=1$, representing the probabilities with which we choose the maps $T_j$. Denote by $m_{\mathbf p}$ the $\mathbf p$-Bernoulli measure on $\Omega^{\mathbb N}$, let $\mu_{\mathbf p}$ be a probability measure on $I$ that is absolutely continuous with respect to the one-dimensional Lebesgue measure $\lambda$ and denote its density by $f_\mathbf p :=\frac{d\mu_\mathbf p}{d\lambda} $. If $\mu_{\mathbf p}$ satisfies for each Borel set $B \subseteq I$ that
\begin{equation}\label{d:invdensity}
\mu_{\mathbf p} (B) = \int_B f_\mathbf p \, d\lambda = \sum_{j \in \Omega} p_j \mu_{\mathbf p}(T_j^{-1}B),
\end{equation}
then the product measure $m_{\mathbf p} \times \mu_{\mathbf p}$ is an invariant probability measure for $R$. Here we call $\mu_{\mathbf p}$ a {\em stationary measure} and $f_{\mathbf p}$ an {\em invariant density} for $R$.

\vskip .2cm
In the literature there exist various sets of conditions under which the existence of such an invariant measure is guaranteed. See for example \cite{Mo,Pe,Bu, GoBo,BaGo1,In}. Here we explicitly mention a special case of the conditions by Inoue from \cite{In} which are simple to state and suit our purposes in the next sections. Let $\Omega \subseteq \mathbb N$, $I \subseteq \mathbb R$ an interval and $\{ T_j:I \to I \}_{j \in \Omega}$ a family of transformations. Let $\mathbf p = (p_j)_{j \in \Omega}$ be a positive probability vector. Assume that the following three conditions hold:
\begin{enumerate}
\item[(a1)] There is a finite or countable interval partition $\{ I_i \}$ of $I$, such that each map $T_j$ is $C^1$ and monotone on the interior of each interval $I_i$.
\end{enumerate}
Let $C$ denote the set of all boundary points of the intervals $I_i$ that are in the interior of $I$. We choose the collection $\{ I_i\}$ as small as possible, so that $C$ contains precisely those points that are a critical point of $T_j$ or $T_j'$ for at least one $j \in \Omega$. We call elements $c \in C$ {\em critical points} for the corresponding random system $R$.
\begin{enumerate}
\item[(a2)] The random system $R$ is {\em expanding on average}, i.e., there exists a constant $0 <\rho < 1$, such that $\sum_{j\in \Omega} \frac{p_j}{|T'_j(x)|} \le \rho$ holds for each $x \in I\setminus C$.
\item[(a3)] For each $j \in \Omega$ and $c\in C$ the map
\[ x \mapsto \begin{cases}
\frac{p_j}{|T_j'(x)|}, & \text{if } x \neq c,\\
0, & \text{otherwise},
\end{cases}\]
is of bounded variation.
\end{enumerate}
It then follows from \cite[Theorem 5.2]{In} that an invariant measure for $R$ of the form $m_{\mathbf p} \times \mu_{\mathbf p}$ with $\mu_{\mathbf p}$ satisfying \eqref{d:invdensity} exists. Let $\mathcal R$ denote the class of random maps $R$ that satisfy these three conditions. We will define random matching for maps in $\mathcal R$, but first we fix some notation on sequences and strings.
 
\vskip .2cm
For each $k > 0$ the set $\Omega^k = \{ \mathbf u = u_1 \cdots u_k \, : \, u_i \in \Omega, \, 1\le i \le k \}$ is the set of all $k$-strings of elements in $\Omega$. We let $\Omega^0 = \{ \epsilon\}$, with $\epsilon$ the empty string. For a finite string $\mathbf u$ let $|\mathbf u|$ denote its length, i.e., $|\mathbf u|=k$ if $\mathbf u \in \Omega^k$. Also, for $1 \le n \le k$ we let $\mathbf u_1^n := u_1 \cdots u_n$ and we set $\mathbf u_1^0 = \epsilon$. Similarly, for an infinite sequence $\omega \in \Omega^{\mathbb N}$ and $n \ge 1$ we use the notation $\omega_1^n := \omega_1 \cdots \omega_n$ with $\omega_1^0= \epsilon$. Finally, we use square brackets to denote cylinder sets, so
\begin{equation}\label{q:cylinder}
[\mathbf u] = \{ \omega \in \Omega^{\mathbb N} \, : \, \omega_1 \cdots \omega_{|\mathbf u|} = \mathbf u \}.
\end{equation}

\vskip .2cm
For $\mathbf u \in \Omega^k$ and $0 \le n \le k$, let
\[ T_{\mathbf u} = T_{u_k} \circ T_{u_{k-1}} \circ \cdots \circ T_{u_1} \quad \text{and} \quad T_{\mathbf u}^n=T_{\mathbf u_1^n} = T_{u_n} \circ T_{u_{n-1}} \circ \cdots \circ T_{u_1}.\]
Note that $T_{\mathbf u}^0=T_{\mathbf u_1^0}=T_{\epsilon}= id$. Similarly if $\omega\in \Omega^\mathbb N$, we let $T^n_{\omega}= T_{\omega_1^n}=T_{\omega_n} \circ T_{\omega_{n-1}} \circ \cdots \circ T_{\omega_1}$ for any $n\ge 0$. For $\mathbf u \in \Omega^k$ the left and right random orbits of the critical points $c \in C$ are
\[ T_{\mathbf u}(c^-) = \lim_{x \uparrow c} T_{\mathbf u}(x) \quad \text{ and } \quad T_{\mathbf u}(c^+) = \lim_{x \downarrow c} T_{\mathbf u}(x).\]
The one-sided derivatives along $\mathbf u$ are given by
\[ T_{\mathbf u}'(c^-) = \lim_{x \uparrow c} \prod_{n=1}^k T_{u_n}'(T_{\mathbf u_1^{n-1}}(x)) \quad \text{ and }\quad T_{\mathbf u}'(c^+) = \lim_{x \downarrow c} \prod_{n=1}^k T_{u_n}'(T_{\mathbf u_1^{n-1}}(x)).\]
We use the abbreviation $p_{\mathbf u} := p_{u_1} \cdots p_{u_k}$ with $p_\epsilon=1$.

\begin{defn}\label{d:rmatching1}{(Random matching)}
A random map $R \in \mathcal R$ has {\em random matching} if for every $c \in C$ there exists an $M=M_c \in \mathbb{N}$ and a set
\[ Y=Y_c \subseteq \Big\{ T_{\omega}^k (c^-) \, : \, \omega \in \Omega^\mathbb N, \, 1 \le k \le M \Big\} \cap \Big\{ T_{\omega}^k (c^+) \, : \, \omega \in \Omega^\mathbb N, \, 1 \le k \le M \Big\} \]
such that for every $\omega \in \Omega^\mathbb N$ there exist $k=k_c(\omega), \ell=\ell_c(\omega) \le M$ with $T^k_\omega (c^-), T^\ell_\omega (c^+)\in Y$.\end{defn}

The main difference with one-dimensional matching as in \eqref{q:onedmatching1} and \eqref{q:onedmatching2} is that in a random system $R$ the critical points have many different random orbits. Definition~\ref{d:rmatching1} states that any random orbit of the left or the right limit of any critical point $c$ passes through the set $Y_c$ at the latest at time $M$. The indices $k, \ell$ are introduced to cater for the possibility that these orbits pass through the set $Y_c$ at different moments. Since all points in $Y_c$ are in the orbit of both $c^-$ and $c^+$, this implies that all random orbits of the left limit meet with some random orbit of the right limit and vice versa. This corresponds to the statement in \eqref{q:onedmatching1}. Note that we do not ask $T^k_\omega (c^-) = T^\ell_\omega (c^+)$.

\begin{defn}\label{d:rmatching2}{(Strong random matching)}
A random map $R \in \mathcal R$ has {\em strong random matching} if it has random matching and if for each $c \in C$ and $y \in Y_c$ the following holds. Set
\[ \begin{split}
\Omega(y)^- =\ & \Big\{ \mathbf u \in \bigcup_{k=1}^M\Omega^k \, : \, \exists\,  \omega \in \Omega^{\mathbb N} \, \text{ with } \mathbf u = \omega_1\cdots \omega_{k_c(\omega)}\,   \text{ and } \, T_{\mathbf u} (c^-) =y \Big\} ,\\
 \Omega(y)^+ =\ & \Big\{ \mathbf u \in \bigcup_{k=1}^M\Omega^k \, : \, \exists\,  \omega \in \Omega^{\mathbb N} \, \text{ with } \mathbf u = \omega_1\cdots \omega_{\ell_c(\omega)}\,   \text{ and } \,  T_{\mathbf u} (c^+)=y \Big\}.
\end{split}\]
Then,
\begin{equation}\label{q:randomderivatives}
\sum_{\mathbf u \in \Omega(y)^-} \frac{ p_{\mathbf u}}{T'_{\mathbf u}(c^-)} = \sum_{\mathbf u \in \Omega(y)^+} \frac{p_{\mathbf u}}{T'_{\mathbf u}(c^+)}.
\end{equation}
\end{defn}

Definition~\ref{d:rmatching2} guarantees that one can choose the times $k, \ell$ such that at those times orbits enter the set $Y$ with the same weighted derivative. This is comparable to \eqref{q:onedmatching2}. Note that
\[ \bigcup_{y \in Y_c} \bigcup_{\mathbf u \in \Omega(y)^-} [\mathbf u] = \Omega^{\mathbb N} = \bigcup_{y \in Y_c} \bigcup_{\mathbf u \in \Omega(y)^+} [\mathbf u],\]
where $[\mathbf u]$ is a cylinder as defined in \eqref{q:cylinder}, so we have indeed captured all random orbits of $c$. Note that Definition~\ref{d:rmatching2} depends on the choices of $k_c(\omega)$ and $\ell_c(\omega)$ for each $c$ in Definition~\ref{d:rmatching1}.

\vskip .2cm
If $\Omega$ consists of one element only, then the random map is actually a deterministic map. In this case Definition~\ref{d:rmatching1} and Definition~\ref{d:rmatching2} reduce to the definitions of one-dimensional matching and strong matching given in \eqref{q:onedmatching1} and \eqref{q:onedmatching2}, so the random definitions extend the deterministic ones.
\vskip .2cm

\subsection{Two examples of families of dynamical systems with random matching}
Below there are two examples of families of random interval maps depending on one parameter. We show that for each of these families there exist parameter intervals such that the systems have strong random matching for every parameter within these intervals. Moreover, within such an interval matching happens in a comparable way, i.e., with the same $M$ and similar sets $Y$. As in the deterministic case, we call these intervals {\em matching intervals}. To ease the notation we use the symbol $\star$ to indicate the set of strings obtained by replacing $\star$ with any $j \in \Omega$. E.g., if $\Omega = \{0,1,2\}$, then $0 \star = \{00, 01, 02\}$. 

\begin{ex}\label{x:randomcf}
{\rm For $\alpha \in (0,1)$ let $T_{\alpha,0}, T_{\alpha,1}: [\alpha-1, \alpha] \rightarrow [\alpha-1, \alpha]$ be the Nakada and Ito-Tanaka $\alpha$-continued fraction transformations, introduced in \cite{Nak81} and in \cite{IT} respectively, which are given by
\[T_{\alpha,0} (x)= \frac{1}{|x|} - \bigg\lfloor \frac{1}{|x|}+1-\alpha \bigg\rfloor \quad \text{and} \quad T_{\alpha,1}(x) = \frac{1}{x} - \bigg\lfloor \frac{1}{x}+1-\alpha \bigg\rfloor ,\]
for $x \neq 0$ and $T_{\alpha,0} (0)=0=T_{\alpha,1} (0)$. The graphs are shown in Figure~\ref{f:randomcf}.
\begin{figure}[h]
\centering
\subfigure[$T_{\alpha,0}$]{\begin{tikzpicture}[scale=3]
\filldraw[ballblue] (-.11494,-.3) rectangle (.11494,.7);

\draw[line width=.4mm, ballblue, smooth, samples =20, domain=-.3:-.27] plot(\x, { -1/ \x -3});
\draw[line width=.4mm, ballblue, smooth, samples =20, domain=-.27:-.2125] plot(\x, { -1/ \x -4});
\draw[line width=.4mm, ballblue, smooth, samples =20, domain=-.2125:-.1756] plot(\x, { -1/ \x -5});
\draw[line width=.4mm, ballblue, smooth, samples =20, domain=-.1756:-.1492] plot(\x, { -1/ \x -6});
\draw[line width=.4mm, ballblue, smooth, samples =20, domain=-.14925:-.12982] plot(\x, { -1/ \x -7});
\draw[line width=.4mm, ballblue, smooth, samples =20, domain=-.12982:-.11494] plot(\x, { -1/ \x -8});
\draw[line width=.4mm, ballblue, smooth, samples =20, domain=.11494:.13001] plot(\x, { 1 / \x -8});
\draw[line width=.4mm, ballblue, smooth, samples =20, domain=.12987:.1493] plot(\x, { 1 / \x -7});
\draw[line width=.4mm, ballblue, smooth, samples =20, domain=.14925:.1756] plot(\x, { 1 / \x -6});
\draw[line width=.4mm, ballblue, smooth, samples =20, domain=.1756:.2125] plot(\x, { 1 / \x -5});
\draw[line width=.4mm, ballblue, smooth, samples =20, domain=.2125:.27] plot(\x, { 1/ \x -4});
\draw[line width=.4mm, ballblue, smooth, samples =20, domain=.27:.37] plot(\x, { 1/ \x -3});
\draw[line width=.4mm, ballblue, smooth, samples =20, domain=.37:.588] plot(\x, { 1/ \x -2});
\draw[line width=.4mm, ballblue, smooth, samples =20, domain=.588:.7] plot(\x, { 1/ \x -1});

\draw(-.305,0)node[left]{\tiny $\alpha-1$}--(-.05,0)node[below]{\tiny $0$}--(.705,0)node[right]{\tiny $\alpha$}(0,-.305)node[below]{\tiny $\alpha-1$}--(0,.705)node[above]{\tiny $\alpha$};
\draw(-.3,-.3)--(-.3,.7)--(.7,.7)--(.7,-.3)--(-.3,-.3);
\draw[dotted] (.588,-.3)node[below]{\small $\frac1{\alpha+1}$}--(.588,.7);
\draw[dotted] (-.27,-.3)--(-.27,.7);
\node[below] at (-.29,-.3) {\small $-\frac1{\alpha+3}$};
\end{tikzpicture}}
\hspace{1cm}
\subfigure[$T_{\alpha,1}$]{\begin{tikzpicture}[scale=3]
\filldraw[blush] (-.1075,-.3) rectangle (.11494,.7);
\draw[line width=.4mm, blush, smooth, samples =20, domain=-.3:-.2326] plot(\x, { 1/ \x +4});
\draw[line width=.4mm, blush, smooth, samples =20, domain=-.2326:-.1887] plot(\x, { 1/ \x +5});
\draw[line width=.4mm, blush, smooth, samples =20, domain=-.1887:-.1587] plot(\x, { 1/ \x +6});
\draw[line width=.4mm, blush, smooth, samples =20, domain=-.1587:-.1367] plot(\x, { 1/ \x +7});
\draw[line width=.4mm, blush, smooth, samples =20, domain=-.137:-.1204] plot(\x, { 1/ \x +8});
\draw[line width=.4mm, blush, smooth, samples =20, domain=-.12048:-.1074] plot(\x, { 1/ \x +9});
\draw[line width=.4mm, blush, smooth, samples =20, domain=.11494:.13001] plot(\x, { 1 / \x -8});
\draw[line width=.4mm, blush, smooth, samples =20, domain=.12987:.1493] plot(\x, { 1 / \x -7});
\draw[line width=.4mm, blush, smooth, samples =20, domain=.14925:.1756] plot(\x, { 1 / \x -6});
\draw[line width=.4mm, blush, smooth, samples =20, domain=.1756:.2125] plot(\x, { 1 / \x -5});
\draw[line width=.4mm, blush, smooth, samples =20, domain=.2125:.27] plot(\x, { 1/ \x -4});
\draw[line width=.4mm, blush, smooth, samples =20, domain=.27:.37] plot(\x, { 1/ \x -3});
\draw[line width=.4mm, blush, smooth, samples =20, domain=.37:.588] plot(\x, { 1/ \x -2});
\draw[line width=.4mm, blush, smooth, samples =20, domain=.588:.7] plot(\x, { 1/ \x -1});

\draw(-.305,0)node[left]{\tiny $\alpha-1$}--(-.05,0)node[below]{\tiny $0$}--(.705,0)node[right]{\tiny $\alpha$}(0,-.305)node[below]{\tiny $\alpha-1$}--(0,.705)node[above]{\tiny $\alpha$};
\draw[dotted] (.588,-.3)node[below]{\small $\frac1{\alpha+1}$}--(.588,.7);
\draw[dotted] (-.2326,-.3)node[below]{\small $\frac1{\alpha-5}$}--(-.2326,.7);
\draw(-.3,-.3)--(-.3,.7)--(.7,.7)--(.7,-.3)--(-.3,-.3);\end{tikzpicture}}
\caption{The Nakada $\alpha$-continued fraction map $T_{\alpha,0}$ in (a) and the Ito-Tanaka $\alpha$-continued fraction map $T_{\alpha,1}$ in (b) for $\alpha = \frac{7}{10} \in \big( \frac{5-\sqrt 13}{2}, \frac{\sqrt 2}{2} \big)$.}
\label{f:randomcf}
\end{figure}
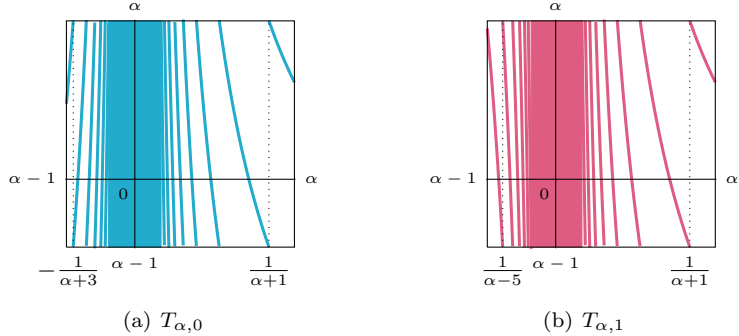

\vskip .2cm
Let $R_\alpha$ denote the corresponding pseudo-skew product on $\{0,1\}^\mathbb N \times [\alpha-1, \alpha]$ and let $\mathbf{p}=(p_0,p_1)$ be a positive probability vector. For $x \in [0,\alpha]$, the two maps coincide and
\[T_{\alpha,0}(x)= T_{\alpha,1}(x)= \frac{1}{x}- n \quad  \text{for } \ x \in \bigg( \frac{1}{\alpha+n}, \frac{1}{\alpha+n-1} \bigg]. \]
For $x \in [\alpha-1,0)$, we have
\[ \begin{array}{ll}
\displaystyle T_{\alpha,0}(x)=  - \frac{1}{x}- n & \text{for } \ \displaystyle x \in \Big[ -\frac{1}{\alpha+n-1}, -\frac{1}{\alpha+n} \Big),\\
\displaystyle T_{\alpha,1}(x)= \frac{1}{x}+ n & \text{for } \ \displaystyle  x \in \Big[ \frac{1}{\alpha-n}, \frac{1}{\alpha-(n+1)} \Big).
\end{array}\]
\vskip .2cm

We first show that for any $\alpha \in \big( \frac{\sqrt{10}-2}{2}, 2-\sqrt 2 \big)$ the map $R_\alpha$ has random matching. For this note that the critical points $c$ are all in the set $\{\frac{1}{\alpha+n}, - \frac{1}{\alpha+n}, \frac{1}{\alpha-n} \, : \, n \in \mathbb{N} \}$. For any positive critical point $c > 0$ and any $j \in \{0,1\}$, $T_j(c^-), T_j(c^+) \in \{\alpha-1, \alpha\}$. For $c < 0$, $c$ is either a critical point for $T_0$ and a continuity point for $T_1$, or a critical point for $T_1$ and a continuity point for $T_0$. Specifically, since $\alpha > \frac12$, for $c= - \frac{1}{\alpha+n}$ we have
\[T_0(c^-)= \alpha, \quad T_0(c^+)= \alpha-1, \quad \text{and} \quad T_1(c^-)=T_1(c^+)=1-\alpha,\]
and for $c= \frac{1}{\alpha-n}$
\[T_1(c^-)= \alpha-1, \quad T_1(c^+)= \alpha, \quad \text{and} \quad T_0(c^-)=T_0(c^+)=1-\alpha.\] 
As a consequence, to show that $R_{\alpha}$ has random matching we only need to consider the orbits of $\alpha-1$ and $\alpha$. Due to the choice of endpoints of the parameter interval $\big( \frac{\sqrt{10}-2}{2}, 2-\sqrt 2 \big)$, the first three orbit points of $\alpha$ and $\alpha-1$ are easily determined. They are given in Figure~\ref{f:randomorbitsCF1}. Hence, if we take $M=3$ and
\[Y= \bigg\{ \frac{5\alpha-3}{1-2\alpha}, \frac{4-7\alpha}{1-2\alpha} \bigg\}, \quad \text{if } c>0\]
and
\[Y= \bigg\{ \frac{5\alpha-3}{1-2\alpha}, \frac{4-7\alpha}{1-2\alpha}, 1-\alpha \bigg\}, \quad \text{if } c<0,\]
then $R_\alpha$ has random matching according to Definition~\ref{d:rmatching1}. 

\begin{figure}[h]
\centering
\resizebox{.85\textwidth}{!}{
\begin{tikzpicture}[->,>=stealth',shorten >=1pt,auto,node distance=2cm,semithick]
\tikzstyle{brightcerulean}=[rectangle,fill=none,draw=brightcerulean,text=black]
\tikzstyle{bluebell}=[rectangle,fill=none,draw=bluebell,text=black]
\tikzstyle{babyblueeyes}=[rectangle,fill=none,draw=babyblueeyes,text=black]
\tikzstyle{blush}=[rectangle,fill=none,draw=blush,text=black]
\node[](A){\footnotesize {\color{white}1} $\alpha$  {\color{white}1}};
\node[](B)[right=0.9cm of A]{\footnotesize $\frac{1-2\alpha}{\alpha}$};  
\node[brightcerulean](C)[right=0.7cm of B]{\footnotesize $\frac{5\alpha-3}{1-2\alpha}$};
\node[blush](K)[below=.7cm of B]{\footnotesize $\frac{4-7\alpha}{1-2\alpha}$};

\path (A) edge node{}(B)
        (B) edge node{\footnotesize 0}(C)
        (B) edge node{\footnotesize 1}(K);
\node[](E)[right=1.1cm of C]{\footnotesize $\alpha-1$};
\node[](F)[right=0.9cm of E]{\footnotesize $\frac{1-2\alpha}{\alpha-1}$};
\node[brightcerulean](G)[right=0.79cm of F]{\footnotesize $\frac{5\alpha -3}{1-2\alpha}$};
\node[](H)[below=0.7cm of E]{\footnotesize $\frac{2\alpha-1}{\alpha-1}$};  
\node[brightcerulean](I)[right=0.7cm of H]{\footnotesize $\frac{5\alpha-3}{1-2\alpha}$};  
\node[blush](J)[below=0.7cm of H]{\footnotesize $\frac{4 - 7\alpha}{1-2\alpha}$};  

\path (E) edge node{\footnotesize 0}(F)
        (F) edge node{}(G)
        (E) edge node{\footnotesize 1}(H)
        (H) edge node{\footnotesize 0}(I)
        (H) edge node{\footnotesize 1}(J);
\end{tikzpicture}
}
\caption{The first three elements in the orbits of $\alpha$ and $ \alpha -1$ under the random continued fraction map $R_{\alpha}$ for $\alpha \in \big( \frac{\sqrt{10}-2}{2}, 2-\sqrt 2 \big)$. The digits above the arrows indicate which one of the maps $T_{\alpha,0}$ or $T_{\alpha,1}$ is applied. If there is no digit, then both maps yield the same orbit point. Orbit points in boxes with the same colour are equal.}
\label{f:randomorbitsCF1}
\end{figure}
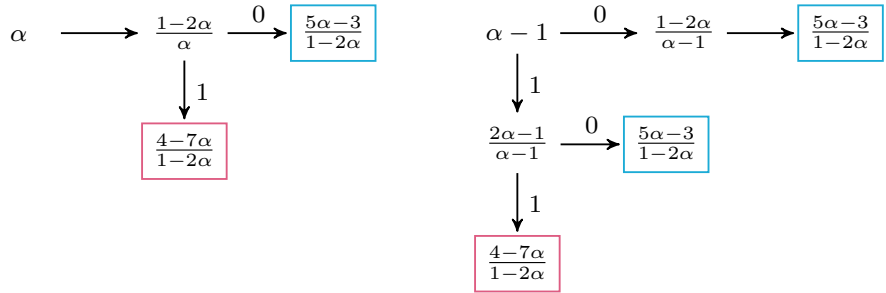

$R_\alpha$ does not satisfy strong random matching with this choice of $Y$. To see this, note that $T_{\alpha,1}'(x) = -\frac1{x^2}$ for all $x$ where the derivative exists, while $T_{\alpha,0}'(x) = -\frac1{x^2}$ if $x>0$ and $T_{\alpha,0}'(x) = \frac1{x^2}$ if $x <0$. Now take for example $c = \frac1{\alpha+n}>0$ and  $y= \frac{4-7\alpha}{1-2\alpha}$. Then
$ \Omega(y)^- = \star 11 = \{ 011, 111\}$ and $\Omega(y)^+ = \star \star 1$. For the quantities from \eqref{q:randomderivatives}, we obtain
\[ \sum_{\mathbf u \in \Omega(y)^-} \frac{p_{\mathbf u}}{T'_\mathbf u}(c^-) = -p_1^2 c^2(2\alpha-1)^2 \quad \text{and} \quad
\sum_{\mathbf u \in \Omega(y)^+} \frac{p_{\mathbf u}}{T'_\mathbf u}(c^+) = -p_1c^2(2\alpha-1)^2,\]
which are not equal for any $p_1 \in (0,1)$.

\vskip .2cm
We now identify a countable number of parameter intervals on which the maps $R_\alpha$ have strong matching with the same exponent $M=4$, i.e., we identify a countable number of {\em matching intervals} for the family $R_\alpha$. For $n \geq 4$ let the interval $J_n:= (\ell_n, r_n )$ be defined by the left and right endpoints
\begin{equation}\label{q:jn}
\ell_n=\frac{n+1-\sqrt{n^2-2n+5}}{2} \quad \text{ and } \quad r_n=\sqrt{\frac{n-2}{n}},
\end{equation} 
respectively. Set $g:=\frac{\sqrt 5-1}{2}$ for the small golden mean and note that $g < \ell_n < r_n$ for all $n \geq 4$ and that $\lim_{n \rightarrow \infty} \ell_n =\lim_{n \rightarrow \infty} r_n = 1$. See Figure \ref{f:intervalsIn} for an illustration of the location of these intervals.
\begin{figure}[h!]
\centering
\includegraphics[trim=200 0 0 0,scale=0.24]{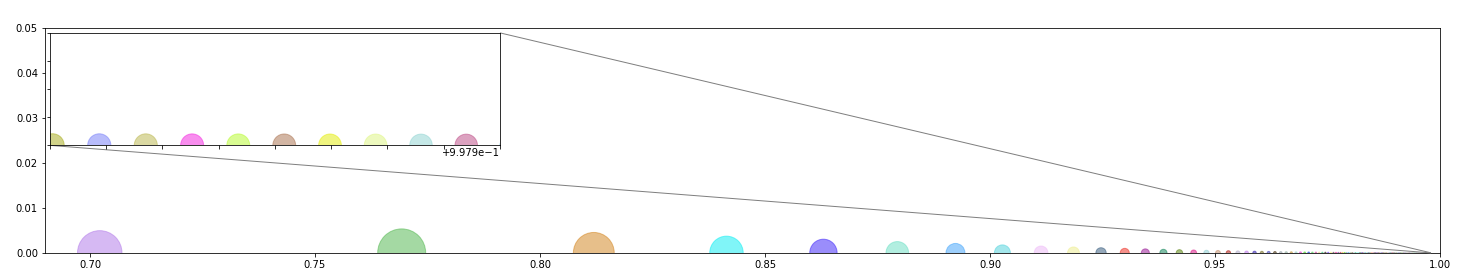}
\caption{The semicircles indicate the locations of the intervals $J_n$.}
\label{f:intervalsIn}
\end{figure}

\vskip .2cm
The intervals $J_n$ are chosen in such a way that we can determine the first three orbit points of $\alpha$ and $\alpha-1$. Let $n \ge 4$ and $\alpha \in J_n$. In particular $\alpha > g$ and for $j=0,1$,
\[T_{\alpha,j}(\alpha)= \frac{1-\alpha}{\alpha}> 0.\]
The point $\ell_n$ is chosen so that
$ \alpha -1 \in (\frac{1}{\alpha-n} ,\frac{\alpha+1}{1-n(\alpha+1)} ) \subseteq (-\frac{1}{\alpha+n-2}, -\frac{1}{\alpha+n-1} )$. Since $\frac{\alpha+1}{1-n(\alpha+1)} < \frac1{\alpha-n-1}$ we get
\[ T_{\alpha,1}(\alpha-1) = \frac{n(\alpha-1)+1}{\alpha-1} \quad \text{ and } \quad T_{\alpha,0} (\alpha-1) = \frac{\alpha(n-1)+2-n}{1-\alpha}.\]
It also implies $ \frac{1-\alpha}{\alpha} \in ( \frac{1}{\alpha+n-2}, \frac{1}{\alpha+n-3} )$. As a consequence, for $l =0,1$, 
\[T_{\alpha, jl} (\alpha) = \frac{\alpha(n-1)+2-n}{1-\alpha} = T_{\alpha,0} (\alpha-1) > 0.\]
We further divide the interval $J_n$. For $k \in \{2, 3, \ldots, n\}$, let 
\[ i_{n,k} = \frac{-4+2n-kn+k + \sqrt{k^2 n^2 -2 k^2 n + k^2 + 4}}{2(n-1)},\]
and note that $J_n \subseteq \cup_{k=2}^{n-1} (i_{n,k+1},i_{n,k}]$. Therefore, for each $\alpha \in J_{n}$ there exists a $k \in \{2, 3, \ldots, n-1\}$ such that $\alpha \in (i_{n,k+1},i_{n,k}]$. The last condition is equivalent to
\begin{equation}\label{e:conditionalpha}
\frac{1}{\alpha+k} < \frac{\alpha(n-1)+2-n}{1-\alpha} \le \frac{1}{\alpha+k-1},
\end{equation}
so that for $\mathbf u  \in \Omega^3$ it holds that 
\[T_{\alpha, \mathbf u}(\alpha) = \frac{1-2k+kn - \alpha (kn-k+1)}{\alpha(n-1)+2-n}.\]
On the other hand, the choice of $r_n$ guarantees that $T_{\alpha,1} (\alpha-1)=\frac{1+n(\alpha-1)}{n} > \frac{1}{\alpha+1}$. Then for $j =0,1$,
\[T_{\alpha,1j} (\alpha-1)=\frac{\alpha(n-1)+2-n}{-1-n(\alpha-1)}.\]
Equation \eqref{e:conditionalpha} holds if and only if 
\[\frac{1}{\alpha+k-1} < \frac{\alpha(n-1)+2-n}{-1-n(\alpha-1)} \leq \frac{1}{\alpha+k-2} \]
is satisfied. In this case, for $l=0,1$
\[T_{\alpha,1jl} (\alpha-1)=\frac{1-2k+kn - \alpha (kn-k+1)}{\alpha(n-1)+2-n}.\]
Figure~\ref{f:randomorbitsCF} shows all the relevant orbit points of $\alpha$ and $\alpha-1$.

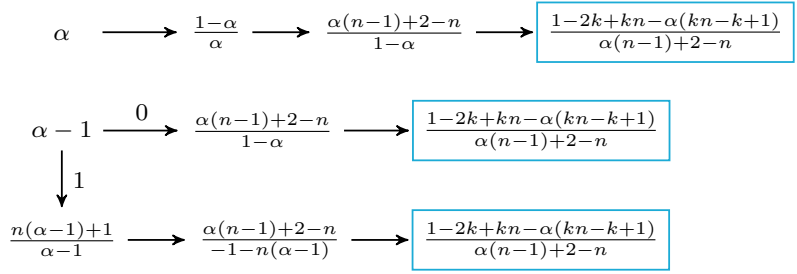
\begin{figure}[h]
\centering

\resizebox{.75\textwidth}{!}{
\begin{tikzpicture}[->,>=stealth',shorten >=1pt,auto,node distance=2cm,semithick]
\tikzstyle{brightcerulean}=[rectangle,fill=none,draw=brightcerulean,text=black]
\tikzstyle{bluebell}=[rectangle,fill=none,draw=bluebell,text=black]
\tikzstyle{babyblueeyes}=[rectangle,fill=none,draw=babyblueeyes,text=black]
\tikzstyle{blush}=[rectangle,fill=none,draw=blush,text=black]
\node[](A){\footnotesize {\color{white}1} $\alpha$  {\color{white}1}};
\node[](B)[right=0.9cm of A]{\footnotesize $\frac{1-\alpha}{\alpha}$};  
\node[](C)[right=0.7cm of B]{\footnotesize $\frac{\alpha(n-1)+2-n}{1-\alpha}$};
\node[brightcerulean](D)[right=0.7cm of C]{\footnotesize $\frac{1-2k+kn - \alpha (kn-k+1)}{\alpha(n-1)+2-n}$};

\path (A) edge node{}(B)
        (B) edge node{}(C)
        (C) edge node{}(D);
\node[](E)[below=.7cm of A]{\footnotesize $\alpha-1$};
\node[](F)[right=0.9cm of E]{\footnotesize $\frac{\alpha(n-1)+2-n}{1-\alpha}$};
\node[brightcerulean](G)[right=0.79cm of F]{\footnotesize $\frac{1-2k+kn - \alpha (kn-k+1)}{\alpha(n-1)+2-n}$};
\node[](H)[below=0.7cm of E]{\footnotesize $\frac{n(\alpha-1)+1}{\alpha-1}$};  
\node[](I)[right=0.7cm of H]{\footnotesize $\frac{\alpha(n-1)+2-n}{-1-n(\alpha-1)}$};  
\node[brightcerulean](J)[right=0.7cm of I]{\footnotesize $\frac{1-2k+kn - \alpha (kn-k+1)}{\alpha(n-1)+2-n}$};  

\path (E) edge node{\footnotesize 0}(F)
        (F) edge node{}(G)
        (E) edge node{\footnotesize 1}(H)
        (H) edge node{}(I)
        (I) edge node{}(J);
\end{tikzpicture}
}
\caption{The first few points in the orbits of $\alpha$ and $ \alpha -1$ under the random continued fraction map $R_{\alpha}$ for $
\alpha \in J_n \cap (i_{n,k+1},i_{n,k}]$.}
\label{f:randomorbitsCF}
\end{figure}

\vskip .2cm
Definition~\ref{d:rmatching1} holds for $\alpha \in J_n \cap (i_{n,k+1},i_{n,k}]$ with $M=4$ and 
\[Y=\bigg\{ \frac{1-2k+kn - \alpha (kn-k+1)}{\alpha(n-1)+2-n }\bigg\}\] 
for any critical point $c>0$. For $c<0$ we add the point $1-\alpha$ to $Y$. Here the values $k_c(\omega)$ and $\ell_c(\omega)$ either equal 1, 3 or 4 according to the number of orbit points in the paths in Figure~\ref{f:randomorbitsCF}. For Definition~\ref{d:rmatching2}, for $c>0$ and $y \in Y$ we have $\Omega(y)^- = \star 0 \star \cup \star 1 \star \star$ and $\Omega(y)^+= \star \star \star \star$, so that
\begin{equation}\nonumber
\begin{split}
\sum_{\mathbf u \in \Omega(y)^-} \frac{p_{\mathbf u}}{T'_{\mathbf u}(c^-)} =\ & (-c^2) p_0 (\alpha-1)^2 \cdot \frac{(\alpha(n-1)+2-n)^2}{-(\alpha-1)^2} \\
& + (-c^2)p_1(-(\alpha-1)^2) \cdot \frac{(1+n(\alpha-1))^2}{-(\alpha-1)^2} \cdot \frac{(\alpha(n-1)+2-n)^2}{-(1+n(\alpha-1))^2} \\ 
=\ & c^2 (\alpha(n-1)+2-n )^2.
\end{split}
\end{equation}
and
\[\sum_{\mathbf u  \in \Omega(y)^+} \frac{ p_{\mathbf u}}{T'_{\mathbf u}(c^+)} = (-c^2) (-\alpha^2) \frac{(1-\alpha)^2}{-\alpha^2} \cdot \frac{(\alpha(n-1)+2-n )^2}{-(1-\alpha)^2} = c^2 (\alpha(n-1)+2-n )^2,\]
implying that also condition \eqref{q:randomderivatives} holds. For $c= -1/(\alpha+n)$ we get $\Omega(1-\alpha)^- =\Omega(1-\alpha)^+ = \{1\}$, $\Omega(y)^- = 0 \star \star \star$ and $\Omega(y)^+ = 0 0 \star \cup \, 0 1 \star \star$, and for $c= 1/(\alpha-n)$ we obtain $\Omega(1-\alpha)^- =\Omega(1-\alpha)^+ = \{0\}$, $\Omega(y)^- =1 0 \star \cup 1 1 \star \star$ and $\Omega(y)^+ = 1 \star \star \star$. In both cases the result follows in a similar fashion. So, the random continued fraction system $R_\alpha$ has strong random matching for any $\mathbf p$ and any $\alpha \in J_n$.
}\end{ex}

\vskip .2cm
Note that in this example the orbits of $\alpha$ meet with some of the orbits of $\alpha-1$ already after two time steps in the point $\frac{\alpha(n-1)+2-n}{1-\alpha}$. Hence,
\[ \frac{\alpha(n-1)+2-n}{1-\alpha} \in \Big\{ T_{\omega}^k (c^-) \, : \, \omega \in \Omega^\mathbb N, \, k \le M \Big\} \cap \Big\{ T_{\omega}^k (c^+) \, : \, \omega \in \Omega^\mathbb N, \, k \le M \Big\}.\]
Therefore, for a critical point $c>0$, we could also take $Y=\big\{\frac{\alpha(n-1)+2-n}{1-\alpha}, \frac{1-2k+kn - \alpha (kn-k+1)}{\alpha(n-1)+2-n } \big\}$ and split the random orbits of $\alpha$ for example in the following way:
 \[  \Omega \bigg( \frac{\alpha(n-1)+2-n}{1-\alpha} \bigg)^+ = \star \star 0 \quad \text{ and } \quad  \Omega \bigg( \frac{1-2k+kn - \alpha (kn-k+1)}{\alpha(n-1)+2-n } \bigg)^+  = \star \star 1 \star.\]
For the orbits passing through $\alpha-1$ we have
 \[ \Omega \bigg( \frac{\alpha(n-1)+2-n}{1-\alpha} \bigg)^- = \star 0 \quad \text{ and } \quad  \Omega \bigg( \frac{1-2k+kn - \alpha (kn-k+1)}{\alpha(n-1)+2-n } \bigg)^- = \star 1 \star \star. \]
One can check that condition \eqref{q:randomderivatives} is satisfied and $R_\alpha$ has strong random matching with this choice of $Y$. Note that in this case many sequences $\omega$ have smaller values $k_c(\omega)$ and $\ell_c(\omega)$ than with $Y=\big\{\frac{1-2k+kn - \alpha (kn-k+1)}{\alpha(n-1)+2-n }\big\}$ and that for some $\omega \in \Omega^\mathbb N$ we do not take $k_c(\omega)$ equal to the first time that the random orbit $T^k_\omega(c^-)$ enters $Y$. For example, for $c >0$ and any $\omega$ with $\omega_3=1$ we have $T^3_\omega(c^+) = \frac{\alpha(n-1)+2-n}{1-\alpha} \in Y$, but we take $k_c(\omega)=4$. The flexibility in the choice of $Y$ and the length of the paths $k_c(\omega)$ and $\ell_c(\omega)$ embedded in Definition~\ref{d:rmatching1} allows one to choose the option that is computationally most convenient.

\begin{ex}\label{x:beta}{\rm Let $\beta = \frac{1+\sqrt 5}{2}$ be the golden mean, so $\beta^2=\beta+1$, and for any $\alpha \in \big( \frac{3\beta-2}{2}, 4\beta-5 \big)$ consider two generalised $\beta$-transformations $T_{\alpha,j}: [-\beta, \beta] \to [-\beta, \beta]$, $j=0,1$, defined by
\[ T_{\alpha,0} (x) = \begin{cases}
\beta x + \alpha, & \text{if } x \in \big[-\beta, - \frac1{\beta} \big),\\
\beta x, & \text{if } x \in \big( -\frac1{\beta}, 1 \big),\\
\beta x -\alpha, & \text{if }x \in (1, \beta],
\end{cases} \quad \text{ and } \quad   T_{\alpha,1} (x) = \begin{cases}
\beta x + \alpha, & \text{if } x \in [-\beta, - 1 ),\\
\beta x, & \text{if } x \in \big( -1, \frac1{\beta} \big),\\
\beta x -\alpha, & \text{if }x \in \big(\frac1{\beta}, \beta\big],
\end{cases}\]
and the maps can be defined however one likes at the discontinuity points. See Figure~\ref{f:beta} for the graphs.
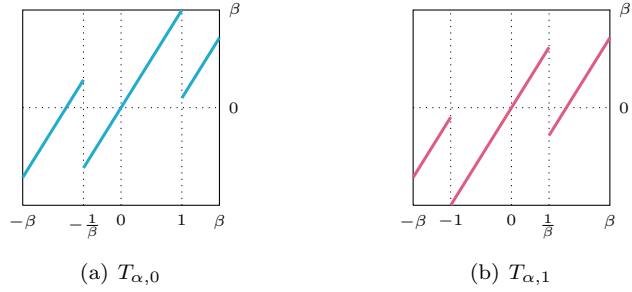
\begin{figure}[h]
\centering
\subfigure[$T_{\alpha,0}$]{
\begin{tikzpicture}[scale=.8]
\draw(-1.618,-1.618)node[below]{\tiny $- \beta$}--(1.618,-1.618)node[below]{\tiny $\beta$}--(1.618,1.618)--(-1.618,1.618)--(-1.618,-1.618);
\draw[dotted] (0,-1.618)--(0,1.618);
\draw[dotted] (1,-1.618)--(1,1.618);
\draw[dotted] (-.618,-1.618)--(-.618,1.618);
\draw[dotted] (-1.618,0)--(1.618,0);

\node[right] at (1.618,0){\tiny $0$};
\node[below] at (-.6,-1.618){\tiny $-\frac{1}{\beta}$};
\node[below] at (1,-1.618){\tiny $1$};
\node[below] at (0,-1.618){\tiny $0$};
\node[right] at (1.618,1.618){\tiny $\beta$};

\draw[line width=.4mm, ballblue, smooth, samples =20, domain=-1.618:-.618] plot(\x, { 1.618* \x +1.46});
\draw[line width=.4mm, ballblue, smooth, samples =20, domain=-.618:1] plot(\x, { 1.618* \x });
\draw[line width=.4mm, ballblue, smooth, samples =20, domain=1:1.618] plot(\x, { 1.618* \x -1.46});
\end{tikzpicture}}
\hspace{1.5cm}
\subfigure[$T_{\alpha,1}$]{
\begin{tikzpicture}[scale=.8]
\draw(-1.618,-1.618)node[below]{\tiny $- \beta$}--(1.618,-1.618)node[below]{\tiny $\beta$}--(1.618,1.618)--(-1.618,1.618)--(-1.618,-1.618);
\draw[dotted] (0,-1.618)--(0,1.618);
\draw[dotted] (.618,-1.618)--(.618,1.618);
\draw[dotted] (-1,-1.618)--(-1,1.618);
\draw[dotted] (-1.618,0)--(1.618,0);

\node[right] at (1.618,0){\tiny $0$};
\node[below] at (-1,-1.618){\tiny $-1$};
\node[below] at (0,-1.618){\tiny $0$};
\node[below] at (.6,-1.618){\tiny $\frac{1}{\beta}$};
\node[right] at (1.618,1.618){\tiny $\beta$};

\draw[line width=.4mm, blush, smooth, samples =20, domain=-1.618:-1] plot(\x, { 1.618* \x +1.46});
\draw[line width=.4mm, blush, smooth, samples =20, domain=-1:.618] plot(\x, { 1.618* \x });
\draw[line width=.4mm, blush, smooth, samples =20, domain=.618:1.618] plot(\x, { 1.618* \x -1.46});
\end{tikzpicture}}
\caption{The maps $T_{\alpha,0}$ and $T_{\alpha,1}$ from Example~\ref{x:beta} for $\alpha \in (\frac{3 \beta-2}{2}, 4 \beta -5 \big)$.}
\label{f:beta}
\end{figure}

\vskip .2cm
Let $R_\alpha$ denote the corresponding random system  and let $\mathbf p = (p_0,p_1)$ be a positive probability vector. Then $C = \big\{ -1, -\frac1{\beta}, \frac1{\beta}, 1 \big\}$. By the symmetry in the maps to show that $R_\alpha$ has matching we only need to consider the points $\frac1{\beta}$ and 1. The parameter interval $\big( \frac{3\beta-2}{2}, 4\beta-5 \big)$ is constructed in such a way that for any $\alpha \in \big( \frac{3\beta-2}{2}, 4\beta-5 \big)$ the initial parts of the random orbits of the left and right limits to $\frac1{\beta}$ and 1 are determined in the following way. For $j=0,1$ and any $\omega \in \{0,1\}^{\mathbb N}$,
\[ \begin{array}{lll}
T_{\alpha,0}(1^-) = \beta, &  T_{\alpha, \omega}(\beta) = \beta^2-\alpha, & T_{\alpha, \omega}^2(\beta) = \beta^2(\beta-\alpha),\\
T_{\alpha,1} (1^-) = T_{\alpha,j}(1^+) = \beta-\alpha, &  T_{\alpha, \omega}(\beta-\alpha) = \beta(\beta-\alpha), & T_{\alpha, \omega}^2(\beta-\alpha) = \beta^2(\beta-\alpha).
\end{array}\]
Hence, for $1 \in C$ we can take $M=k_1(\omega)=\ell_1(\omega)=3$ for each $\omega$, $Y = \{ \beta^2(\beta-\alpha)\}$ and one easily checks the conditions of both Definition~\ref{d:rmatching1} and Definition~\ref{d:rmatching2}.

\vskip .2cm
For $\frac1{\beta}$ the orbits are more complicated. Firstly, $T_{\alpha,j}\big( \frac1{\beta}^-\big)=1 = T_{\alpha,0} \big( \frac1{\beta}^+\big)$ and $T_{\alpha,1}\big( \frac1{\beta}^+\big)=1-\alpha$. We saw the orbit of 1 above, so we concentrate on the orbit of $1-\alpha$. We have $T_{\alpha,j} (1-\alpha) = \beta(1-\alpha) \in \big( -1, -\frac1{\beta} \big)$, so $T_{\alpha,0}(\beta(1-\alpha)) = \beta(\beta-\alpha)$ and $T_{\alpha,1} (\beta(1-\alpha)) = \beta^2(1-\alpha)$. The next couple of iterations are depicted in Figure~\ref{f:randomorbits}, where we have used the property that $\beta^2=\beta+1$ to compute the orbit points.
\begin{figure}[h]
\centering
\resizebox{1\textwidth}{!}{
\begin{tikzpicture}[->,>=stealth',shorten >=1pt,auto,node distance=2.2cm,semithick]
\tikzstyle{brightcerulean}=[rectangle,fill=none,draw=brightcerulean,text=black]
\tikzstyle{bluebell}=[rectangle,fill=none,draw=bluebell,text=black]
\tikzstyle{babyblueeyes}=[rectangle,fill=none,draw=babyblueeyes,text=black]
\tikzstyle{blush}=[rectangle,fill=none,draw=blush,text=black]
\node[](A){\footnotesize $1-\alpha$};
\node[](B)[right=0.7cm of A]{\footnotesize $\beta(1-\alpha)$};
\node[](C)[right=0.7cm of B]{\footnotesize $\beta(\beta-\alpha)$};
\node[blush](D)[right=0.7cm of C]{\footnotesize $\beta^2(\beta-\alpha)$};
\node[](E)[right=0.7cm of D]{\footnotesize $\beta^3(\beta-\alpha)$};
\node[](F)[right=0.7cm of E]{\footnotesize $\beta^4(\beta-\alpha)$};
\node[bluebell](G)[right=0.7cm of F]{\footnotesize $\beta^5(\beta-\alpha)-\alpha$};
\node[](H)[below=0.7cm of B]{\footnotesize $\beta^2(1-\alpha)$};  
\node[](I)[below=0.7cm of H]{\footnotesize $\beta^3(1-\alpha)+\alpha$};  
\node[](J)[below=0.7cm of I]{\footnotesize $\beta^4(1-\alpha)+ \beta \alpha$};  
\node[](M)[below=0.7cm of J]{\footnotesize $\beta^5(1-\alpha)+\beta^2 \alpha$};  
\node[babyblueeyes](N)[below=0.7cm of M]{\footnotesize $\beta^6(1-\alpha)+ \beta^3\alpha +\alpha = \beta^5(\beta-\alpha)-\beta\alpha$};  
\node[](K)[right=0.7cm of J]{\footnotesize $\beta^5(1-\alpha)+ \beta^2 \alpha + \alpha$};  
\node[bluebell](L)[right=0.7cm of K]{\footnotesize $\beta^6(1-\alpha)+ \beta^3 \alpha + \beta \alpha = \beta^5(\beta-\alpha)-\alpha$};  
\node[](O)[below=0.7cm of E]{\footnotesize $\beta^4(\beta-\alpha)-\alpha$};
\node[babyblueeyes](P)[below=0.7cm of O]{\footnotesize $\beta^5(\beta-\alpha)- \beta \alpha$};

\path (A) edge node{}(B)
        (B) edge node{\small 0}(C)
        (C) edge node{}(D)
        (D) edge node{}(E)
        (E) edge node{\small 0}(F)
        (F) edge node{}(G)
        (B) edge node{\small 1}(H)
        (H) edge node{}(I)
        (I) edge node{}(J)
        (J) edge node{\small 1}(M)
        (M) edge node{}(N)
        (J) edge node{\small 0}(K)
        (K) edge node{}(L)
        (E) edge node{\small 1}(O)
        (O) edge node{}(P);
\end{tikzpicture}}
\caption{The first couple of points in the orbit of $1 - \alpha$ under the random generalised $\beta$-transformation from Example~\ref{x:beta}. We have boxed $\beta^2(\beta-\alpha)$, since this point also appears in all random orbits of 1.}
\label{f:randomorbits}
\end{figure}
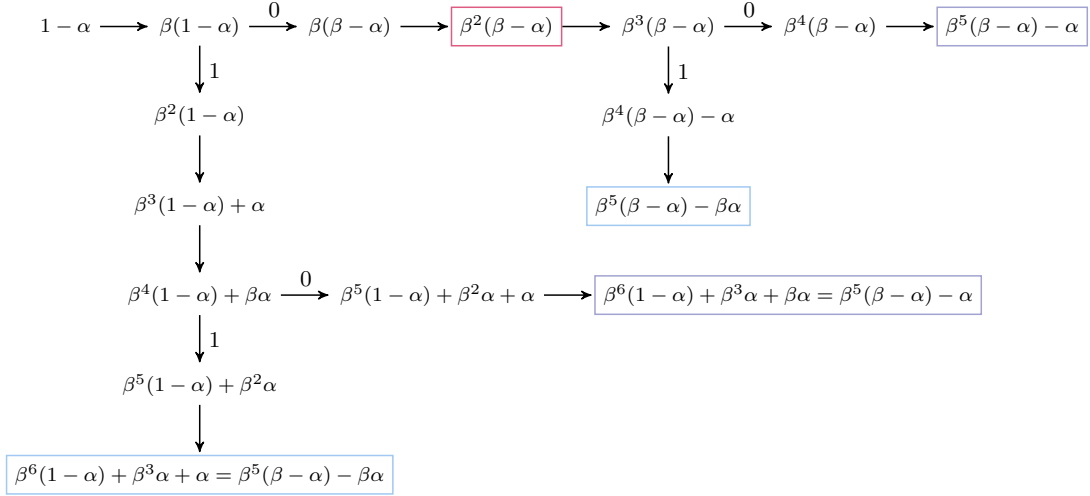

\vskip .2cm
Take $M= k_{\frac1{\beta}}(\omega)=\ell_{\frac1{\beta}}(\omega)= 7$ for each $\omega$ and set $Y = \{ \beta^5(\beta-\alpha)-\alpha,  \beta^5(\beta-\alpha)-\beta\alpha = \beta^6 -3\beta^3\alpha \}$. Then,
\[ \Omega(\beta^5(\beta-\alpha)-\alpha)^+ = 0 \star \star \star \star 0 \star \cup 1 \star 0 \star \star 0 \star \cup 1 \star 1 \star \star 0 \star \]
and $\Omega(\beta^5(\beta-\alpha)-\alpha)^- = \star \star \star \star \star 0 \star$. Hence,
\[  \sum_{\mathbf u  \in \Omega(\beta^5(\beta-\alpha)-\alpha)^+} \frac{ p_{\mathbf u}}{T'_{\mathbf u}(\frac1{\beta}^+)} =\frac{p_0^2 + p_1p_0^2 + p_1^2p_0}{\beta^7} = \frac{p_0}{\beta^7} = \sum_{\mathbf u \in \Omega(\beta^5(\beta-\alpha)-\alpha)^-} \frac{p_{\mathbf u}}{T'_{\mathbf u}(\frac1{\beta}^-)}.\] 
A similar computation gives \eqref{q:randomderivatives} for $\beta^5(\beta-\alpha)-\beta\alpha$, so $R_\alpha$ has strong random matching.
}\end{ex}

\vskip .2cm
Note that also in this example the orbits of $1^+$ meet with some of the orbits of $1^-$ earlier, in this case already after one step. Hence, we could also take $Y_1 = \{ \beta- \alpha, \beta^2(\beta-\alpha) \}$ and split the random orbits as follows:
\[ \Omega(\beta-\alpha)^+ = \{ 1\} = \Omega(\beta-\alpha)^- \quad \text{and} \quad  \Omega(\beta^2(\beta-\alpha))^+ = 0 \star \star = \Omega(\beta^2 (\beta-\alpha))^-.\]
Then for some $\omega$ the values $k_1(\omega), \ell_1(\omega)$ are lower, but we have to check condition \eqref{q:randomderivatives} for two points instead of one. For the critical point $\frac1{\beta}$ we could use $Y = \{ \beta-\alpha, \beta^2 (\beta-\alpha), \beta^5(\beta-\alpha)-\alpha, \beta^5(\beta-\alpha)-\beta\alpha \}$ or also $Y = \{\beta^2 (\beta-\alpha), \beta^5(\beta-\alpha)-\alpha, \beta^5(\beta-\alpha)-\beta\alpha \}$. By the flexibility in the choice of $Y$ given by Definition~\ref{d:rmatching1} one can choose the set $Y$ that is most convenient.
Theorem~\ref{t:pwconstant} below explains the need for condition \eqref{q:randomderivatives} in Definition~\ref{d:rmatching2}.

\subsection{Random matching for piecewise affine systems}
In case each map $T_j: I \to I$ is piecewise affine on a finite partition $c_0 < c_1 < \ldots < c_N$ the conditions (a1) and (a3) are automatically satisfied and under some additional assumptions strong random matching has consequences for invariant densities. For this result we consider a subset of the collection of random maps $\mathcal R$. We define the subset $\mathcal R_A \subset \mathcal R$ to be the set of random systems in $\mathcal R$ that satisfy the following additional assumptions:
\begin{itemize}
\item[(c1)] There exists a finite interval partition $\{ I_i\}_{1 \le i \le N}$ of $I=[c_0, c_N]$ given by the points $c_0 < c_1 < \ldots < c_N$, such that each map $T_j:I \to I$, $j \in \Omega$, is piecewise affine with respect to this partition. In other words, for each $j \in \Omega$ and $1 \le i \le N$ we can write $T_j|_{(c_{i-1},c_i)}(x) = k_{i,j} x + d_{i,j}$ for some constants $k_{i,j}, d_{i,j}$.
\item[(c2)] For each $1 \le i \le N$ there is an $1 \le n \le N$, such that
\begin{equation}\label{q:diagonal}
\frac{\sum_{j \in \Omega} \frac{p_j}{k_{i,j}}d_{i,j}}{1- \sum_{j \in \Omega} \frac{p_j}{k_{i,j}}} \neq \frac{\sum_{j \in \Omega} \frac{p_j}{k_{n,j}}d_{n,j}}{1- \sum_{j \in \Omega} \frac{p_j}{k_{n,j}}}.
\end{equation}
\item[(c3)] For each $1 \le i \le N$, $\sum_{j \in \Omega} \frac{p_j}{k_{i,j}} \neq 0$.
\end{itemize}
Using the results from \cite{KM}, we will show that for $R \in \mathcal R_A$ the following holds.

\begin{thm}\label{t:pwconstant}
Let $R \in \mathcal R_A$. If $R$ has strong random matching, then there exists an invariant probability measure $m_{\mathbf p} \times \mu_{\mathbf p}$ for $R$ with $\mu_{\mathbf p}$ absolutely continuous with respect to Lebesgue and such that its density $f_{\mathbf{p}}$ is piecewise constant. If moreover every map $T_j$ is expanding, i.e., if $|k_{i,j}|>1$ for each $1 \le i \le N$ and $j \in \Omega$, then any invariant probability density $f_{\mathbf p}$ of $R$ is piecewise constant.
\end{thm}

\vskip .2cm
Assumptions (c2) and (c3) are used in \cite{KM} to prove that for systems in $\mathcal R_A$ there exists an invariant probability density function that can be written as an infinite sum of indicator functions. We use this fact in the proof below. These conditions, which are not very restrictive, guarantee that the method from \cite{KM} works, but they might not be necessary for the results from \cite[Theorem 4.1]{KM} and Theorem~\ref{t:pwconstant}. In fact, the deterministic analog of Theorem~\ref{t:pwconstant}, which can be found in \cite[Theorem 1.2]{BCMP}, does not have a condition like \eqref{q:diagonal}. Their proof uses an induced system with a full branched return map instead. One could try to transfer the proof of \cite[Theorem 1.2]{BCMP} to the setting of random interval maps to avoid (c2) and (c3). Then, the recent results from Inoue in \cite{Ino20} on first return time functions for random systems seem relevant. These results show, however, that an induced system for a random interval map will become position dependent instead of i.i.d., which might make such an extension not so straightforward.

\begin{proof}
The set of critical points of $R$ is given by $C = \{ c_1, \ldots, c_{N-1}\}$. Any random map $R \in \mathcal R_A$ satisfies the conditions of \cite[Theorem 4.1]{KM}. Thus, there exists an invariant probability measure $m_{\mathbf p} \times \mu_{\mathbf p}$ for $R$ with a probability density $f_\mathbf p$ for $\mu_{\mathbf p}$ of the form
\begin{equation}\label{q:densityf}
 f_\mathbf p = \sum_{i=1}^{N-1} \gamma_i \sum_{k \ge 1} \sum_{\mathbf u \in \Omega^k} \Big( \frac{p_{\mathbf u}}{T'_{\mathbf u}(c_i^-)} 1_{[c_0, T_{\mathbf u}(c_i^-))}- \frac{p_{\mathbf u}}{T'_{\mathbf u}(c_i^+)} 1_{[c_0, T_{\mathbf u}(c^+_i))}\Big),
 \end{equation}
for some constants $\gamma_i$ depending only on the critical points $c_i$. Fix an $i$ and let $M,Y$ be such that $R$ satisfies the conditions of Definition~\ref{d:rmatching1} and Definition~\ref{d:rmatching2} for $c_i$. Then by \eqref{q:randomderivatives}
\[ \sum_{y \in Y} \Big( \sum_{\mathbf u \in \Omega(y)^-} \frac{p_{\mathbf u}}{T'_{\mathbf u}(c_i^-)} 1_{[c_0, T_{\mathbf u}(c_i^-))} - \sum_{\mathbf u \in \Omega(y)^+} \frac{p_{\mathbf u}}{T'_{\mathbf u}(c_i^+)} 1_{[c_0, T_{\mathbf u}(c_i^+))}  \Big)=0.\]
For each $1 \le i \le N-1$ and each $1 \le k \le M$, let
\[ \Omega^{i,k}_- = \{ \mathbf u \in \Omega^k \, : \, \exists \, \omega \in \Omega^\mathbb N \text{ with } \mathbf u = \omega_1 \cdots \omega_k \text{ and } k < k_{c_i}(\omega) \}\]
and similarly
\[ \Omega^{i,k}_+ = \{ \mathbf u \in \Omega^k \, : \, \exists \, \omega \in \Omega^\mathbb N \text{ with } \mathbf u = \omega_1 \cdots \omega_k \text{ and } k < \ell_{c_i}(\omega) \}.\]
Then $f_\mathbf p$ can be written as
\[  f_\mathbf p = \sum_{i=1}^{N-1} \gamma_i \sum_{k = 1}^M \Big( \sum_{\mathbf u \in \Omega^{i,k}_-} \frac{p_{\mathbf u}}{T'_{\mathbf u}(c_i^-)} 1_{[c_0, T_{\mathbf u}(c_i^-))}-  \sum_{\mathbf u \in \Omega^{i,k}_+} \frac{p_{\mathbf u}}{T'_{\mathbf u}(c_i^+)} 1_{[c_0, T_{\mathbf u}(c^+_i))}\Big).\]
Hence $f_\mathbf p$ is constant on each interval in the finite partition of $I$ specified by the orbit points in the set
\[ \bigcup_{i=1}^{N-1} \bigcup_{k=1}^M \Big(\{ T_{\mathbf u} (c_i^-) \, :  \,  \mathbf u \in \Omega^{i,k}_- \} \cup \{ T_{\mathbf u} (c_i^+) \, : \, \mathbf u \in \Omega^{i,k}_+ \} \Big).\]
This gives the first part of the result.

\vskip .2cm
For the second part, note that under the additional assumption that $|k_{i,j}|>1$ for all $i,j$ the map $R$ satisfies the conditions of \cite[Theorem 5.3]{KM}. As a consequence, any invariant density $f_\mathbf p$ of $R$ can be written as in \eqref{q:densityf} for some values $\gamma_i$. This proves the theorem.
\end{proof}

\begin{ex}{\rm
The random generalised $\beta$-transformations $R_\alpha$ from Example~\ref{x:beta} satisfy all conditions of Theorem~\ref{t:pwconstant}. 
Hence, for any $\alpha \in \big( \frac{3\beta-2}{2}, 4\beta-5 \big)$ any invariant density of the random system $R_\alpha$ is piecewise constant.
}\end{ex}

\section{Random signed binary transformations and expansions}
In the second part of this article we use strong random matching to study the frequency of the digit 0 in the signed binary expansions produced by a family of random system of piecewise affine maps. We first define this family and its relation to binary expansions.

\subsection{The family of random symmetric doubling maps}
\begin{figure}[h]
\centering
\subfigure[$D_a$]{\begin{tikzpicture}[scale=1.2]
\draw(-1,-1)node[below]{\tiny $-1$}--(1,-1)node[below]{\tiny $1$}--(1,1)--(-1,1)--(-1,-1);
\draw[dotted] (0,-1)--(0,1);
\draw[dotted] (.5,-1)--(.5,1);
\draw[dotted] (-.5,-1)--(-.5,1);
\draw[dotted] (-1,0)--(1,0);
\draw[dotted] (-1,-1)--(1,1);
\draw[dashed] (.35,-1)--(.35,1)(-.35,-1)--(-.35,1)(-1,.7)--(-.7,.7)(-1,-.7)--(-.7,-.7);
\draw[red] (-.7,-.7)--(.7,-.7)--(.7,.7)--(-.7,.7)--(-.7,-.7);

\node[right] at (1,0){\tiny $0$};
\node[below] at (-1,-1){\tiny $-1$};
\node[below] at (-.5,-1){\tiny $-\frac{1}{2}$};
\node[below] at (.5,-1){\tiny $\frac{1}{2}$};
\node[below] at (.35,-1){\tiny $a$};
\node[right] at (1,1){\tiny $1$};
\node[left] at (-1,.7){\tiny $2a$};
\node[left] at (-1,-.7){\tiny $-2a$};

\draw[line width=.4mm, ballblue, smooth, samples =20, domain=-1:-.35] plot(\x, { 2* \x +1});
\draw[line width=.4mm, gainsboro, smooth, samples =20, domain=-.35:0] plot(\x, { 2* \x +1});
\draw[line width=.4mm, ballblue, smooth, samples =20, domain=-.35:.35] plot(\x, { 2* \x });
\draw[line width=.4mm, gainsboro, smooth, samples =20, domain=-.5:-.35] plot(\x, { 2* \x });
\draw[line width=.4mm, gainsboro, smooth, samples =20, domain=.35:.5] plot(\x, { 2* \x });
\draw[line width=.4mm, gainsboro, smooth, samples =20, domain=0:.35] plot(\x, { 2* \x -1});
\draw[line width=.4mm, ballblue, smooth, samples =20, domain=.35:1] plot(\x, { 2* \x -1});
\end{tikzpicture}}
\hspace{1.5cm}
\subfigure[$S_{\frac1{2a}}$]{\begin{tikzpicture}[scale=1.2]
\draw(-1,-1)node[below]{\tiny $- 1$}--(1,-1)node[below]{\tiny $1$}--(1,1)--(-1,1)--(-1,-1);
\draw[dotted] (0,-1)--(0,1);
\draw[dashed] (.5,-1)--(.5,1);
\draw[dashed] (-.5,-1)--(-.5,1);
\draw[dotted] (-1,0)--(1,0)(-1,-1)--(1,1);

\node[right] at (1,0){\tiny $0$};
\node[below] at (-1,-1){\tiny $-1$};
\node[below] at (-.5,-1){\tiny $-\frac{1}{2}$};
\node[below] at (1,-1){\tiny $1$};
\node[below] at (0,-1){\tiny $0$};
\node[below] at (.5,-1){\tiny $\frac{1}{2}$};
\node[right] at (1,1){\tiny $1$};
\node[right] at (1,.5){\tiny $\frac12$};

\draw[line width=.4mm, ballblue, smooth, samples =20, domain=-1:-.5] plot(\x, { 2* \x +10/7});
\draw[line width=.4mm, ballblue, smooth, samples =20, domain=-.5:.5] plot(\x, { 2* \x  });
\draw[line width=.4mm, ballblue, smooth, samples =20, domain=.5:1] plot(\x, { 2* \x -10/7});
\end{tikzpicture}}
\caption{The maps $D_a$ and $S_\frac1{2a}$ for $a=\frac{7}{20}$. The grey lines indicate the remainder of the maps $x \mapsto 2x+1$, $x \mapsto 2x$ and $x \mapsto 2x-1$. The red box in (a) shows the attractor of the map $D_c$. The map inside the box in (a) is a rescaled version of the map in (b).}
\label{f:doublingdet}
\end{figure}
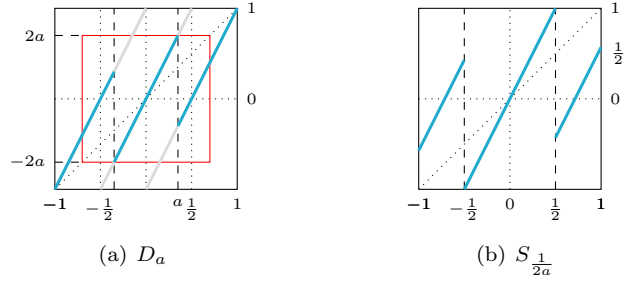

A signed binary expansion of a number $x \in [-1,1]$ can be obtained by iterating any piecewise affine map $D:[-1,1] \to [-1,1]$ that is given by $D(x) = 2x - d$ with $d \in \{-1,0,1\}$ on each of its intervals of monotonicity. One can for example take any $a \in \big[ \frac14, \frac12 \big]$ and then define the symmetric map
\[ D_a(x) = \begin{cases}
2x +1, & \text{if } -1 \le x < -a,\\
2x, & \text{if } -a \le x \le a,\\
2x-1, & \text{if } a < x \le 1.
\end{cases}\]
By setting $d_n(x) = d$, $d \in \{-1,0,1\}$, if $D_a^n(x) = 2D^{n-1}_a(x)-d$, one obtains
\[ x = \frac{d_1(x)}{2} + \frac{D_a(x)}{2} = \cdots = \frac{d_1(x)}{2} + \cdots + \frac{d_n(x)}{2^n} + \frac{D_a^n(x)}{2^n} \to \sum_{n \ge 1} \frac{d_n(x)}{2^n},\]
so this gives a signed binary expansion of $x$. The family of maps $\{D_a\}_{\frac14 \le a \le \frac12}$ is the object of study in \cite{DK}. As can be seen from Figure~\ref{f:doublingdet}(a) the interval $[-2a,2a]$ is an attractor for the dynamics of $D_a$. Since this interval depends on $a$, in \cite{DK} the authors decided to work instead with the measurably isomorphic family $\{S_\alpha\}_{1 \le \alpha \le 2}$ given by
\begin{equation}\label{q:Seta}
S_\alpha(x) = \begin{cases}
2x +\alpha, & \text{if } -1 \le x < -\frac12,\\
2x, & \text{if } -\frac12 \le x \le \frac12,\\
2x-\alpha, & \text{if } \frac12 < x \le 1,
\end{cases}
\end{equation}
see Figure~\ref{f:doublingdet}(b), which transfers the dependence on the parameter from the domain $[-1,1]$ to the branches of the maps.

\vskip .2cm
While each deterministic map produces for each number in its domain a single signed binary expansion, one can define random dynamical systems that produce for Lebesgue almost all numbers uncountably many different signed binary expansions. The family of random maps $\{ R_\alpha\}$, which we define next, extends the family of deterministic maps $\{ S_\alpha \}$. So the dependence on the parameter is visible in the branches of the maps instead of in the domains.

\begin{figure}[h]
\centering
\subfigure[$R_1$]{\begin{tikzpicture}[scale=1.1]
\draw(-1,-1)node[below]{\tiny $-1$}--(1,-1)node[below]{\tiny $1$}--(1,1)--(-1,1)--(-1,-1);
\draw[dotted] (0,-1)--(0,1);
\draw[dotted] (.5,-1)--(.5,1);
\draw[dotted] (-.5,-1)--(-.5,1);
\draw[dotted] (-1,0)--(1,0);
\filldraw[gray!30,nearly transparent] (-.5,-1) rectangle (.5,1);

\node[right] at (1,0){\tiny $0$};
\node[below] at (-1,-1){\tiny $-1$};
\node[below] at (-.5,-1){\tiny $-\frac{1}{2}$};
\node[below] at (.5,-1){\tiny $\frac{1}{2}$};
\node[right] at (1,1){\tiny $1$};

\draw[line width=.4mm, bluebell, smooth, samples =20, domain=-1:-.5] plot(\x, { 2* \x +1});
\draw[line width=.4mm, babyblueeyes, smooth, samples =20, domain=-.5:0] plot(\x, { 2* \x +1});
\draw[line width=.4mm, babyblueeyes, smooth, samples =20, domain=0:.5] plot(\x, { 2* \x });
\draw[line width=.4mm, blush, smooth, samples =20, domain=-.5:0] plot(\x, { 2* \x });
\draw[line width=.4mm, blush, smooth, samples =20, domain=0:.5] plot(\x, { 2* \x -1});
\draw[line width=.4mm, bluebell, smooth, samples =20, domain=.5:1] plot(\x, { 2* \x -1});
\end{tikzpicture}}
\hspace{1cm}
\subfigure[$R_{\frac32}$]{\begin{tikzpicture}[scale=1.1]
\draw(-1,-1)node[below]{\tiny $- 1$}--(1,-1)node[below]{\tiny $1$}--(1,1)--(-1,1)--(-1,-1);
\draw[dotted] (0,-1)--(0,1);
\draw[dotted] (.5,-1)--(.5,1);
\draw[dotted] (.25,-1)--(.25,1);
\draw[dotted] (-.5,-1)--(-.5,1);
\draw[dotted] (-.25,-1)--(-.25,1);
\draw[dotted] (-1,0)--(1,0);
\filldraw[gray!30,nearly transparent] (-.5,-1) rectangle (-.25,1);
\filldraw[gray!30,nearly transparent] (.25,-1) rectangle (.5,1);

\node[right] at (1,0){\tiny $0$};
\node[below] at (-1,-1){\tiny $-1$};
\node[below] at (-.5,-1){\tiny $-\frac{1}{2}$};
\node[below] at (-.25,-1){\tiny $-\frac{1}{4}$};
\node[below] at (1,-1){\tiny $1$};
\node[below] at (0,-1){\tiny $0$};
\node[below] at (.5,-1){\tiny $\frac{1}{2}$};
\node[below] at (.25,-1){\tiny $\frac{1}{4}$};
\node[right] at (1,1){\tiny $1$};
\node[right] at (1,.5){\tiny $\frac12$};

\draw[line width=.4mm, bluebell, smooth, samples =20, domain=-1:-.5] plot(\x, { 2* \x +3/2});
\draw[line width=.4mm, babyblueeyes, smooth, samples =20, domain=-.5:-.25] plot(\x, { 2* \x +3/2});
\draw[line width=.4mm, bluebell, smooth, samples =20, domain=-.25:.25] plot(\x, { 2* \x  });
\draw[line width=.4mm, blush, smooth, samples =20, domain=-.5:-.25] plot(\x, { 2* \x  });
\draw[line width=.4mm, babyblueeyes, smooth, samples =20, domain=.25:.5] plot(\x, { 2* \x  });
\draw[line width=.4mm, blush, smooth, samples =20, domain=.25:.5] plot(\x, { 2* \x -3/2});
\draw[line width=.4mm, bluebell, smooth, samples =20, domain=.5:1] plot(\x, { 2* \x -3/2});
\end{tikzpicture}}
\hspace{1cm}
\subfigure[$R_2$]{\begin{tikzpicture}[scale=1.1]
\draw(-1,-1)node[below]{\tiny $- 1$}--(1,-1)node[below]{\tiny $1$}--(1,1)--(-1,1)--(-1,-1);
\draw[dotted] (0,-1)--(0,1);
\draw[dotted] (.5,-1)--(.5,1);
\draw[dotted] (-.5,-1)--(-.5,1);
\draw[dotted] (-1,0)--(1,0);

\node[right] at (1,0){\tiny $0$};
\node[below] at (-1,-1){\tiny $-1$};
\node[below] at (-.5,-1){\tiny $-\frac{1}{2}$};
\node[below] at (1,-1){\tiny $1$};
\node[below] at (0,-1){\tiny $0$};
\node[below] at (.5,-1){\tiny $\frac{1}{2}$};
\node[right] at (1,1){\tiny $1$};

\draw[line width=.4mm, bluebell, smooth, samples =20, domain=-1:-.5] plot(\x, { 2* \x +2});
\draw[line width=.4mm, bluebell, smooth, samples =20, domain=-.5:.5] plot(\x, { 2* \x  });
\draw[line width=.4mm, bluebell, smooth, samples =20, domain=.5:1] plot(\x, { 2* \x -2});
\end{tikzpicture}}
\caption{The maps $T_{\alpha,0}$ and $T_{\alpha,1}$ for $\alpha=1$ in (a), $\alpha=\frac32$ in (b), and $\alpha=2$ in (c). The blue lines correspond to $T_{\alpha, 0}$, the pink ones to $T_{\alpha,1}$ and the purple ones to both.}
\label{f:rsb}
\end{figure}
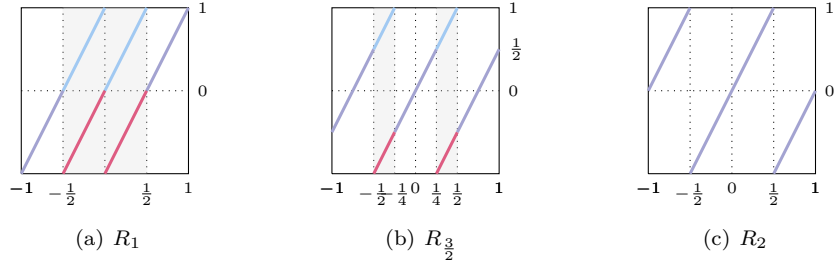

\vskip .2cm
Let $\Omega = \{0,1\}$ and define for $j \in \Omega$ and each parameter $\alpha \in [1,2]$ the maps $T_j=T_{\alpha,j}:[-1,1]\to [-1,1]$ by
\begin{equation}\label{q:t0t1}
T_{\alpha,0}(x) = \begin{cases}
2 x + \alpha, & \text{if } x \in \big[-1, \frac{1-\alpha}{2}\big], \\
2 x, & \text{if } x \in \big(\frac{1-\alpha}{2},\frac12\big],\\
2 x - \alpha, & \text{if } x \in \big(\frac12,1\big],
\end{cases}
\, \text{ and } \, T_{\alpha,1}(x) = \begin{cases}
2 x + \alpha, & \text{if } x \in \big[-1,-\frac12 \big), \\
2 x, & \text{if } x \in \big[-\frac12,\frac{\alpha-1}{2}\big),\\
2 x - \alpha, & \text{if } x \in \big[\frac{\alpha-1}{2},1\big].
\end{cases}
\end{equation}
See Figure \ref{f:rsb} for three examples. The maps $T_{\alpha,0}$ and $T_{\alpha,1}$ differ on the intervals $\big[-\frac12, 1-2\alpha \big]$ and $\big[2\alpha-1, \frac12 \big]$, which are indicated by the grey areas in the pictures. Let $R=R_\alpha : \Omega^\mathbb N \times [-1,1] \to \Omega^\mathbb N \times [-1,1]$ be the random system obtained from $T_{\alpha,0}$ and $T_{\alpha,1}$, i.e., 
\[ R_\alpha (\omega, x) = \big( \sigma (\omega), T_{\alpha, \omega_1} (x) \big),\]
where $\sigma$ is the left shift on $\Omega^\mathbb N$. We call the systems $R_\alpha$ {\em random symmetric doubling maps} and  the subscript $\alpha$ will sometimes be suppressed if it does not lead to confusion.

\vskip .2cm
Fix an $\alpha \in [1,2]$. Recall from \eqref{q:cylinder} that we use square brackets to denote the cylinder sets in $\Omega^\mathbb N$. Let $\pi: \Omega^\mathbb N \times [-1,1] \to [-1,1]$ denote the canonical projection onto the second coordinate and set
\begin{equation}\label{q:randomdigits}
s_n(\omega,x) = \begin{cases}
-1, & \text{if } R^{n-1} (\omega,x) \in \Omega^\mathbb N \times \big[-1, -\frac12 \big) \cup [0] \times \big[ -\frac12, \frac{1-\alpha}{2} \big],\\
0, & \text{if } R^{n-1} (\omega,x) \in [1] \times  \big[ -\frac12, \frac{1-\alpha}{2} \big]\\
& \hspace{2.5cm} \cup\, \Omega^\mathbb N \times \big(\frac{1-\alpha}{2}, \frac{\alpha-1}{2}\big) \cup [0] \times \big[\frac{\alpha-1}{2}, \frac12\big],\\
1, & \text{if } R^{n-1}(\omega,x) \in [1]\times \big[ \frac{\alpha-1}{2}, \frac12 \big] \cup \Omega^\mathbb N \times \big( \frac12,1\big].
\end{cases}
\end{equation}
Then
\[ \pi (R^n(\omega,x)) = 2 \pi(R^{n-1}(\omega,x)) - s_n(\omega,x)\alpha,\]
so that just as in the deterministic case by iteration we obtain
\[ x = \frac{s_1(\omega,x)\alpha}{2} + \cdots + \frac{s_n(\omega,x)\alpha}{2^n} + \frac{\pi (R^n(\omega,x))}{2^n} \to \alpha \sum_{n \ge 1} \frac{s_n(\omega,x)}{2^n}.\]
In other words, iterations of the random system $R$ give a signed binary expansion for the pair $(\omega,x)$.

\vskip .2cm
Note that for each $x \in [-1,1]$ there is an $\omega \in \Omega^\mathbb N$, such that $\pi(R_\alpha^n(\omega, x)) = S^n_\alpha (x)$, where $S_\alpha$ is the map in the family $\{ S_\alpha\}$ from \cite{DK}. In particular, the random signed binary expansions produced by the family $\{R_\alpha\}$ include, among many others, the SSB expansions. The randomness of the system allows us to choose (up to a certain degree) where and when we want to have a digit 0. Below we investigate the frequency of the digit 0 in typical expansions produced by the maps $R$. We do so by applying Birkhoff's Ergodic Theorem for invariant measures for $R$ of the form $m \times \mu$ with $m$ a Bernoulli measure and $\mu$ absolutely continuous with respect to the Lebesgue measure. For that we need to investigate the density of such measures $\mu$.

\subsection{Prevalent matching for random symmetric doubling maps}
For any $\alpha \in (1,2]$ the common partition on which $T_0$ and $T_1$ are monotone is given by the points
\[c_0=-1, \quad c_1=- \frac{1}{2} , \quad c_2=\frac{1-\alpha}{2}, \quad c_3=\frac{\alpha-1}{2}, \quad c_4= \frac12, \quad c_5=1.\]
Set, in accordance with (a1), 
\[I_1= [c_0,c_1), \quad I_2= [c_1,c_2], \quad I_3= (c_2,c_3), \quad I_4= [c_3,c_4], \quad I_5= (c_4,c_5],\]
then $C= \{c_1, c_2, c_3, c_4\}$. For $0 < p < 1$, use $\mathbf p=(p_0, p_1)$ to denote the probability vector with $p_0=p$ and $p_1=1-p$. Since $T_0$ and $T_1$ from \eqref{q:t0t1} are both piecewise affine with slope 2, we have $\frac{p_j}{|T_j'(x)|} = \frac{p_j}{T_j'(x)} = \frac{p_j}{2}$, $j=0,1$. So the random system $R$ satisfies conditions (a1), (a2), (a3), i.e., $R \in \mathcal R$. Due to the symmetry in the map, to verify whether $R$ has strong random matching it is enough to check the conditions of Definitions~\ref{d:rmatching1} and Definitions~\ref{d:rmatching2} for the points $1 = T_0(c_4)^-$ and $1-\alpha = T_0(c_4)^+$.

\vskip .2cm
Before we proceed with a description of the matching behaviour of the family of random systems $\{R_\alpha\}$, we first recall the results from \cite[Propositions 2.1 and 2.3]{DK} on matching for the family of deterministic symmetric doubling maps $\{ S_\alpha\}$, see \eqref{q:Seta}. Let
\begin{equation}\label{q:Meta}
M_\alpha = \inf\Big\{ n \ge 0 \, : \, \frac12 < S_\alpha^n(1) < \alpha - \frac12 \Big\} +1.
\end{equation} 
Then according to \cite[Propositions 2.1 and 2.3]{DK} for all $\alpha \in [1,2]$,
\begin{equation}\label{q:distancealpha}
S_\alpha^k (1-\alpha)= S_\alpha^k (1) -\alpha \quad \text{for } k < M_\alpha
\end{equation}
and for Lebesgue almost all $\alpha \in [1,2]$ in fact $M_\alpha < \infty$ and
\begin{equation}\nonumber
S_\alpha^{M_\alpha+1} \left(\frac12^-\right) = S_\alpha^{M_\alpha} (1) = S_\alpha^{M_\alpha} (1-\alpha) = S_\alpha^{M_\alpha+1} \left(\frac12^+\right).
\end{equation}
In other words, for Lebesgue almost all parameters $\alpha \in [1,2]$ the map $S_\alpha$ has matching with matching exponent $M=M_\alpha+1$ that is determined by the first time the orbit of 1 enters the interval $\big( \frac12, \alpha -\frac12 \big)$. Moreover, $S_\alpha^{M_\alpha-1} (1-\alpha) < - \frac12$ for all $\alpha$, $M_{\alpha}=1$ for $\alpha \in \big[ \frac32, 2\big]$ and $M_{\alpha}>1$ for $\alpha \in \big(1, \frac32\big)$. Due to the constant slope and the same matching exponent $M_{\alpha}$ of the left and right limits, in this case matching implies strong matching.

\begin{nrem}
The discrepancy between $M_\alpha+1$ here and $M_\alpha$ as matching exponent in \cite{DK} comes from the fact that in \cite{DK} the orbits are considered as starting from $1$ and $1-\alpha$, whereas in \eqref{q:onedmatching1} and \eqref{q:onedmatching2} we followed the convention in \cite{BCMP} and start at the critical point $c = \frac12$ instead.
\end{nrem}

From this we deduce the following small lemma.
\begin{lem}\label{l:i4i5}
For $\alpha \in [1, 2]$ and for all $k < M_\alpha-1$, either $S^k(1) \in I_4$ and $S^k(1)-\alpha \in I_1$ or $S^k(1) \in I_5$ and $S^k(1)-\alpha \in I_2$.
\end{lem}

\begin{proof}
From \eqref{q:distancealpha} it follows for all $k < M_\alpha-1$ that $S^k(1) -\alpha \ge -1$, implying that
\[\alpha-1 \leq S^k(1) \leq 1 \quad \text{and} \quad -1 \leq S^k(1) -\alpha \leq 1-\alpha.\]
The fact that $S^k(1) \in I_4 \cup I_5$ follows since $\frac{\alpha-1}{2} \le \alpha-1$. If $S^k(1) \in I_4$, then $S^k(1)-\alpha \le \frac12-\alpha < -\frac12$, so $S^k(1)-\alpha \in I_1$. Suppose $S^k(1) \in I_5$. If $S^k(1) -\alpha < -\frac12$, this would imply that $S^k(1) \in \big( \frac12, \alpha-\frac12\big)$, contradicting the definition of $M_\alpha$ in \eqref{q:Meta}. Hence, $S^k(1) -\alpha \in I_2$.
\end{proof}

The next result states that a random equivalent of \eqref{q:distancealpha} holds for $\alpha \in \big(1, \frac32 \big)$.

\begin{prop}\label{p:distancealpha}
For all $\alpha \in [1,2]$, $0 \leq k \leq M_\alpha$ and $\mathbf u \in \Omega^k$, it holds that $T_{\mathbf u}(1),T_{\mathbf u} (1-\alpha) \in \{ S^k(1), S^k(1)-\alpha \}$.
\end{prop}

\begin{proof}
First consider $\alpha \in \big[\frac32,2\big]$. Then $M_\alpha =1$ and the result trivially holds. Fix an $\alpha \in \big(1,\frac32\big)$. Since $T_0$ and $T_1$ agree on $I_5$ we can find a sequence $\widehat{\omega}\in \Omega^{\mathbb N}$ with $\widehat{\omega}_1=0$ that gives
\[T_{\widehat{\omega}_1^k} (1)= S^k(1) \quad \text{ for all } k \ge 0.\]
Note that $1 \in I_5$ and from $\alpha \in \big(1, \frac32 \big)$ we get $1-\alpha \in I_2$, so
\begin{equation}\label{e:fit1}
T_0 (1-\alpha) =T_0(1)=T_1(1) =2-\alpha = S(1) \quad \text{and} \quad T_1 (1-\alpha)= 2-2\alpha = S(1-\alpha).
\end{equation}
Hence, from the first iterate on, the orbits of $1$ and $1-\alpha$ under the deterministic map $S$ are contained in the orbit of $1-\alpha$ under the random map $R$. To prove the statement, we therefore only have to consider $T_\omega^n (1-\alpha)$ for any $\omega \in \Omega^{\mathbb N}$ and $n \ge 1$. In particular \eqref{e:fit1} implies that 
\[T_{\widehat{\omega}_1^k} (1-\alpha)=S^k(1)\] 
for all $k \geq 1$. We prove the statement by induction.

\vskip .1cm
The statement obviously holds for $k=0$ and by \eqref{e:fit1} also for $k=1$. Let $1\le n < M_\alpha$ and suppose the statement holds for all $k \le n$. Then
\[T_{\widehat{\omega}_1^n} (1-\alpha)= S^n(1) \quad \text{and} \quad T_{\omega_1^n} (1-\alpha) \in \{S^n(1), S^n(1)-\alpha\} \quad \text{for all } \omega \in \Omega^\mathbb N.\]
By Lemma~\ref{l:i4i5} there are three cases.\\
1. If $S^n(1) \in I_4$, then $S^{n+1} (1) = 2S^n(1)$ and $S^n(1)-\alpha \in I_1$. So for the random images we get
\[ T_0(S^n(1))= 2S^n(1), \quad T_1(S^n(1))=2S^n(1)-\alpha,\]
and
\[ T_0(S^n(1)-\alpha)= T_1(S^n(1)-\alpha)=2S^n(1)-\alpha.\]
2. If $S^n(1) \in I_5$ and $S^n(1)-\alpha \in I_2$, then
\[T_0 (S^n(1))= T_1(S^n(1)) =2S^n(1) - \alpha= S^{n+1}(1)\]
and
\[T_0 (S^n(1)-\alpha) = 2S^n(1)-\alpha, \quad T_1 (S^n(1)) = 2S^n(1)-2\alpha.\]
3. If $S^n (1) \in I_5$ and $ S^n(1) \in I_1$ (so $n=M_\alpha-1$), then 
\[T_0(S^n(1)) = T_1(S^n(1)) = 2S^n-\alpha = S^{n+1}(1)\]
and
\[ T_0(S^n(1)-\alpha) = T_1(S^n(1)-\alpha) = 2S^n(1)-\alpha = S^{n+1}(1).\]
Hence, for all $\mathbf u \in \Omega^n$ and $j=0,1$, $T_{\mathbf u j}(1-\alpha) \in \{ S^{n+1}(1), S^{n+1}(1)-\alpha\}$, which gives the result.
\end{proof}

From this proposition we can deduce that matching is prevalent for the family $\{ R_\alpha\}$ and we can find the precise matching times. We first prove the following lemma, stating that all the orbit points $S^n(1), S^n(1-\alpha)$ up to the moment of matching are different.
\begin{lem}\label{l:difforbit}
For each $k < M_\alpha$ the set $\{ S^n(1), S^n(1-\alpha) \, : \, 0 \le n \le k\}$ has $2(k+1)$ elements. 
\end{lem}

\begin{proof}
Since $k < M_\alpha$ it follows from \eqref{q:distancealpha} that $S^n(1) \neq S^n(1-\alpha)$ for each $n$. It also cannot hold that there are $0 \le n < k<M_\alpha$ such that $S^n(1) = S^k(1-\alpha)$ or $S^k(1) = S^n(1-\alpha)$, since this would imply that $|S^k(1)-S^n(1)|= \alpha$ and that would contradict the fact that $S^n(1), S^k(1) \in I_4 \cup I_5$. This leaves the possibility that there are $0 \le n < k<M_\alpha$ such that $S^n(1) = S^k(1)$, i.e., that the orbit of 1 under $S$ is ultimately periodic, or $S^n(1-\alpha) = S^k(1-\alpha)$. Assume $S^n(1) = S^k(1)$ for some $n < k$. It follows that $S^n(1-\alpha) = S^n(1)-\alpha = S^k(1)-\alpha = S^k(1-\alpha)$, so the orbit of $1-\alpha$ is also ultimately periodic and by Proposition~\ref{p:distancealpha} all these orbit points lie at distance $\alpha$ of the corresponding orbit points of 1. This contradicts the fact that $\alpha$ is a matching parameter. Hence, the set $\{ S^n(1), S^n(1-\alpha) \, : \, 0 \le n \le k\}$ has $2(k+1)$ elements.
\end{proof}

\begin{thm}\label{t:mae}
For Lebesgue almost all parameters $\alpha \in [1,2]$ the map $R_{\alpha}$ has strong random matching with $M=M_{\alpha}+1$, where $M_\alpha$ is given by \eqref{q:Meta}, and $Y = \{ S^{M_\alpha}(1)\}$. Moreover, $R_\alpha$ does not satisfy the conditions of strong random matching for any $K < M$.
\end{thm}

\begin{proof}
First consider $\alpha \in \big[\frac32,2\big]$. Then $T_j (1-\alpha) = 2-\alpha =  T_j(1)$ for $j=0,1$, so random matching occurs for $R$ with $M=2$ and $Y=\{ 2-\alpha\}$ and both parts of the theorem hold.

\vskip .2cm
Now, fix $\alpha \in [1, \frac32)$ such that $S=S_\alpha$ has matching. Then, $S^k(1) \neq S^k(1-\alpha)$ for $1 \leq k < M_\alpha$ and $S^{M_\alpha-1}(1) \in \big(\frac12, \alpha-\frac12\big)$, so that $S^{M_\alpha}(1)=2S^{M_\alpha-1}(1)-\alpha$. By Proposition~\ref{p:distancealpha} for each $\mathbf u \in \Omega^{M_\alpha-1}$ either
\[ T_{\mathbf u}(1-\alpha) = S^{M_\alpha-1}(1) > \frac12\]
or 
\[ T_{\mathbf u} (1-\alpha) = S^{M_\alpha-1}(1) - \alpha <- \frac12.\]
In both cases this leads to $T_{\mathbf u j}(1-\alpha) = 2S^{M_\alpha-1}(1) - \alpha$ for both $j=0,1$. The same statement holds for $T_{\mathbf u}(1)$, so that for $c = \frac12$ we therefore have
\[ T_{1\mathbf u j} \left( \frac12^- \right) = T_{0\mathbf u j} \left( \frac12^+ \right) = T_{1\mathbf u j} \left( \frac12^+ \right) = T_{\mathbf u j} (1-\alpha) = S^{M_\alpha}(1)\]
and
\[ T_{0\mathbf u j} \left( \frac12^- \right) = T_{\mathbf u j} (1) = S^{M_\alpha}(1).\]
Hence, we can take $Y_{\frac12} = \{ S^{M_\alpha} (1)\}$. Since this set contains one element only and the maps $T_j$ have the same constant slope, condition \eqref{q:randomderivatives} from Definition~\ref{d:rmatching2} follows immediately. The first part of the theorem now follows since the deterministic maps $S_\alpha$ have matching for Lebesgue almost all parameters $\alpha$. For the critical points $c \neq \frac12$ the statement follows by symmetry.

\vskip .2cm
For the second part we assume for $\alpha \in  [1, \frac32)$ that $S=S_\alpha$ has matching and we proceed by contradiction. Therefore, assume that $R_\alpha$ satisfies the conditions of Definition~\ref{d:rmatching1} and Definition~\ref{d:rmatching2} for $c=\frac12$ for some minimal $1 \le K <M=M_\alpha+1$. Suppose that $S^n(1) \in Y_{\frac12}$ for some $n < K-1$. By Lemma~\ref{l:difforbit} any $\mathbf u$ for which $T_\mathbf u \big(\frac12^\pm \big) = S^n(1)$ has length $|\mathbf u|=n+1$. Together with \eqref{q:randomderivatives} and the fact that the maps $T_{\alpha,0}$ and $T_{\alpha,1}$ both have constant slope 2, this implies that
\begin{equation}\label{q:omega+-}
\sum_{\mathbf u \in \Omega(S^n(1))^-} p_{\mathbf u} = \sum_{\mathbf u \in \Omega(S^n(1))^+} p_{\mathbf u}.
\end{equation}
For any $\mathbf u \in \Omega^{n+1} \setminus \Omega(S^n(1))^-, \mathbf u' \in \Omega^{n+1} \setminus \Omega(S^n(1))^+$ we have by Proposition~\ref{p:distancealpha} that $T_\mathbf u \big(\frac12^-\big) = T_{\mathbf u'}\big(\frac12^+\big) = S^n(1)-\alpha$. Furthermore,
\[\begin{split}
1 = \sum_{\mathbf u \in \Omega^{n+1}} p_{\mathbf u} =\ & \sum_{\mathbf u \in  \Omega(S^n(1))^-} p_{\mathbf u} + \sum_{\mathbf u \in \Omega^{n+1} \setminus \Omega(S^n(1))^-} p_{\mathbf u}\\
=\ & \sum_{\mathbf u \in  \Omega(S^n(1))^+} p_{\mathbf u} + \sum_{\mathbf u \in \Omega^{n+1} \setminus \Omega(S^n(1))^+} p_{\mathbf u}.
\end{split}\]
From \eqref{q:omega+-} and Proposition~\ref{p:distancealpha} we see that
\[  \sum_{\mathbf u \in \Omega(S^n(1-\alpha))^-} p_{\mathbf u} = \sum_{\mathbf u \in \Omega^{n+1} \setminus \Omega(S^n(1))^-} p_{\mathbf u} = \sum_{\mathbf u \in \Omega^{n+1} \setminus \Omega(S^n(1))^+} p_{\mathbf u} = \sum_{\mathbf u \in \Omega(S^n(1-\alpha))^+} p_{\mathbf u}.\]
This implies that the conditions of Definition~\ref{d:rmatching1} and Definition~\ref{d:rmatching2} hold with $M_{\frac12} = n+1$ and $Y_{\frac12} = \{ S^n(1), S^n(1-\alpha)\}$, contradicting the minimality of $K$. In a similar way we can exclude the possibility that $S^n(1-\alpha) \in Y$ for $n < K-1$. Since there is an $\widetilde{\omega}\in \Omega^{\mathbb N}$ such that for each $k <M-1$, $T_{\widetilde{\omega}_1^k} (1-\alpha) = S^k (1-\alpha) = S^k(1)-\alpha$, it must hold that 
\[ Y_{\frac12} = \{ S^{K-1}(1), S^{K-1}(1)-\alpha\}.\]

\vskip .2cm
To conclude the proof we show that for this set $Y_{\frac12}$ condition \eqref{q:randomderivatives} cannot hold. By the constant slope, condition \eqref{q:randomderivatives} can be rephrased as
\begin{equation}\label{e:dercon}
\begin{dcases}
\sum_{\mathbf u \in \Omega(S^{K-1}(1))^-}  p_{\mathbf u} - \sum_{\mathbf u \in \Omega(S^{K-1}(1))^+}  p_{\mathbf u}&=0, \\
\sum_{\mathbf u \in \Omega(S^{K-1}(1)-\alpha)^-}  p_{\mathbf u} - \sum_{\mathbf u \in \Omega(S^{K-1}(1)-\alpha)^+}  p_{\mathbf u}&=0.
\end{dcases}
\end{equation}
and by Lemma~\ref{l:difforbit} any $\mathbf u \in \Omega(S^{K-1}(1))^\pm \cup \Omega(S^{K-1}(1)-\alpha)^\pm$ has length $K$. Since $K<M_{\alpha}+1$, so $K-2 < M_\alpha-1$, Lemma~\ref{l:i4i5} tells us that there are only two possibilities:
\begin{itemize}
\item[1.] $S^{K-2}(1) \in I_4$ and $S^{K-2}(1)-\alpha \in I_1$;
\item[2.] $S^{K-2}(1) \in I_5$ and $S^{K-2}(1)-\alpha \in I_2$.
\end{itemize}
If case 1.~holds, then $T_0(S^{K-2}(1)) = S^{K-1}(1)$ and
\[ T_1(S^{K-2}(1))= T_0(S^{K-2}(1)-\alpha)=T_1(S^{K-2}(1)-\alpha)=S^{K-1}(1)-\alpha,\]
so that \eqref{e:dercon} becomes
\begin{equation}\nonumber
\begin{dcases}
\sum_{\mathbf u \in  \Omega(S^{K-2}(1))^-}  p_{\mathbf u}p_0 - \sum_{\mathbf u \in  \Omega(S^{K-2}(1))^+}  p_{\mathbf u} p_0 =0, \\
\sum_{\mathbf u \in \Omega(S^{K-2}(1))^-}  p_{\mathbf u}p_1 + \sum_{\mathbf u \in (\Omega(S^{K-2}(1)-\alpha)^-}  p_{\mathbf u} - \sum_{\mathbf u \in \Omega(S^{K-2}(1))^+}  p_{\mathbf u}p_1 - \sum_{\mathbf u \in (\Omega(S^{K-2}(1)-\alpha)^+}  p_{\mathbf u} =0.
\end{dcases}
\end{equation}
The last system of equations implies
\begin{equation}\nonumber
\begin{dcases}
\sum_{\mathbf u \in \Omega(S^{K-2}(1))^-}  p_{\mathbf u}- \sum_{\mathbf u \in \Omega(S^{K-2}(1))^+}  p_{\mathbf u} &=0, \\
\sum_{\mathbf u \in (\Omega(S^{K-2}(1)-\alpha)^-}  p_{\mathbf u}  - \sum_{\mathbf u \in (\Omega(S^{K-2}(1)-\alpha)^+}  p_{\mathbf u}&=0, 
\end{dcases}
\end{equation}
which contradicts the minimality of $K$. For the second case, the same contradiction is obtained in a similar way.
\end{proof}

\begin{nrem}
From the previous result we see that matching occurs for the random systems $R_\alpha$ for the same parameters $\alpha$ and at the same time as for the deterministic systems $S_\alpha$. \cite{DK} contains a complete description of the matching intervals of the maps $S_\alpha$. The interval $[1,2]$ can be divided into intervals of parameters for which matching of the maps $S_\alpha$ occurs after the same number of steps. By the above, these matching intervals also apply to the systems $R_\alpha$.
\end{nrem}

\subsection{An expression for the invariant density}
Let $\lambda$ be the Lebesgue measure on $[-1,1]$. The existence of an invariant measure of the form $m_{\mathbf p} \times \mu_{\mathbf p}$ with $\mu_{\mathbf p} \ll \lambda$ for the random symmetric doubling maps $R_\alpha$ is guaranteed by the results of \cite{Pe, Mo}. Furthermore, since $T_0$ is expanding and has a unique absolutely continuous invariant measure, it follows from \cite[Corollary 7]{Pe} that also for $R_\alpha$ there is a unique measure $m_{\mathbf p} \times \mu_{\mathbf p}$ and that $R_\alpha$ is ergodic with respect to this measure. To show that $R_\alpha \in \mathcal R_A$, we check conditions (c1), (c2), (c3). (c1) is immediate and (c3) follows from the constant slope 2 of the maps $T_{\alpha,0}$ and $T_{\alpha,1}$. We check condition~\eqref{q:diagonal}. Note that for any $\alpha \neq 2$,
\[ \frac{\sum_{j \in \Omega} \frac{p_j}{k_{3,j}}d_{3,j}}{1- \sum_{j \in \Omega} \frac{p_j}{k_{3,j}}} = \frac{\sum_{j \in \Omega} \frac{p_j}{2}0}{1- \sum_{j \in \Omega} \frac{p_j}{2}}=0\]
and
\[ \frac{\sum_{j \in \Omega} \frac{p_j}{k_{1,j}}d_{1,j}}{1- \sum_{j \in \Omega} \frac{p_j}{k_{1,j}}} = \frac{\sum_{j \in \Omega} \frac{p_j}{2}\alpha}{1- \sum_{j \in \Omega} \frac{p_j}{2}} =2 \alpha\neq 0.\]
Then \cite[Theorem 5.3]{KM} implies that an explicit formula for the density of this measure can be found via the algebraic procedure in \cite{KM} and from Theorem~\ref{t:pwconstant} and Theorem~\ref{t:mae} we know that for Lebesgue almost all parameters $\alpha$ this density is piecewise constant. We will execute the procedure from \cite{KM} and start by introducing the same notation as in \cite{KM}. Since $\Omega$ consists of two elements only, from now on we will just use $p$ as an index instead of $\mathbf p$ whenever appropriate. 

\vskip .2cm
Denote by $a_{i,j}$ and $b_{i,j}$ the left and right limits at each critical point $c_i \in C$, i.e., for $1 \leq i \leq 4$ and $j \in \Omega$:
\[a_{i,j}= T_j(c_i^-)= \lim_{x \uparrow c_i} T_j(x), \qquad \text{and} \qquad b_{i,j}= T_j(c_i^+)= \lim_{x \downarrow c_i} T_j(x).\]
The images of the critical points are then given by 
\[ \begin{array}{llll}
a_{1,0}=a_{1,1}=b_{1,0}=\alpha-1,  & \hspace{.6cm} b_{1,1}=-1,  &  \hspace{.6cm} a_{2,0}=1,   & \hspace{.6cm} a_{2,1}=b_{2,0}=b_{2,1}=1-\alpha,\\
a_{3,0}=a_{3,1}=b_{3,0}=\alpha-1,  & \hspace{.6cm} b_{3,1}=-1,  & \hspace{.6cm}  a_{4,0}=1,  &  \hspace{.6cm} a_{4,1}=b_{4,0}=b_{4,1}=1 - \alpha.
\end{array}\]

\vskip .2cm
\noindent For $y \in [-1,1]$ and $1 \leq n \leq 4$ set
\begin{equation}\label{q:kin}
\KI_n(y)= \sum_{k \geq 1 } \sum_{\mathbf u\in \Omega^k}  \frac{p_{\mathbf u}}{2^k} 1_{I_n} (T_{\mathbf u^{k-1}_1}(y)).
\end{equation}
The quantity $\KI_n(y)$ weighs the number of times the random orbits of $y$ enters the interval $I_n$. The weight depends on the length and probability of each path $\omega \in \Omega^\mathbb{N}$ leading the point $y$ to $I_n$. The {\em fundamental matrix} $A= (A_{n,i})$ of $R$ is the $5 \times 4$ matrix with entries
\begin{equation}\nonumber
A_{n,i}=  \begin{cases}
\vspace{.2cm}
\displaystyle \, \sum_{j\in \Omega} p_j ( 1 + \KI_n(a_{i,j})- \KI_n(b_{i,j}) ), &\text{for } n = i, \\
\vspace{.2cm}
\displaystyle \, \sum_{j\in \Omega} p_j (\KI_n(a_{i,j})-\KI_n(b_{i,j}) -1) , &\text{for } n = i+1, \\
\displaystyle \, \sum_{j\in \Omega} p_j (\KI_n(a_{i,j})- \KI_n(b_{i,j})), &\text{else.} 
 \end{cases}
 \end{equation}
Since for $R$ there is a unique invariant probability measure $m_p \times \mu_p$ with $\mu_p \ll \lambda$, \cite[Theorem 5.3]{KM} implies that the null space of the matrix $A$ is one-dimensional. According to \cite[Theorem 4.1]{KM} there is a unique vector $\gamma = (\gamma_1, \gamma_2, \gamma_3, \gamma_4) \in \mathbb R^4 \setminus \{ \mathbf 0 \}$ with $A \gamma =\mathbf 0$ and such that the probability density $f_p$ of $\mu_p$ has the form \eqref{q:densityf}. Using the values of $a_{i,j}$ and $b_{i,j}$ computed above, we can reduce this to
\begin{equation}\label{e:fdensity}
\begin{split}
f_p =  (\gamma_1 & + \gamma_3) \sum_{k \ge 0} \sum_{\mathbf u \in \Omega^k} \frac{p_{1\mathbf u}}{2^{k+1}} \big( 1_{[-1, T_\mathbf u(\alpha-1))} - 1_{[-1,T_\mathbf u(-1))} \big)\\
& + (\gamma_2 + \gamma_4) \sum_{k \ge 0} \sum_{\mathbf u \in \Omega^k} \frac{p_{0\mathbf u}}{2^{k+1}} \big( 1_{[-1, T_\mathbf u(1))} - 1_{[-1,T_\mathbf u(1-\alpha))} \big).
\end{split}
\end{equation}
By symmetry 
to determine $f_p$ it is enough to know the random orbits of $1$ and $1-\alpha$ only. From \eqref{e:fdensity} we see that the density is piecewise constant when the orbits of $1$ and $1-\alpha$ are finite or when they merge with the same weight. In the former case the map admits a Markov partition, the latter case happens if $R$ exhibits strong random matching. We focus on the second situation, since we know from Theorem~\ref{t:mae} that this holds for Lebesgue almost all parameters.

\vskip .2cm
Fix an $\alpha \in [1,2]$ such that $R$ presents strong random matching. Let $M$ be as in Theorem~\ref{t:mae}. Then for each $i,j,n$,
\[ \KI_n(a_{i,j}) - \KI_n(b_{i,j}) = \sum_{k = 1}^{M-1} \sum_{\mathbf u \in \Omega^k} \frac{p_\mathbf u}{2^k} \big(1_{I_n} (T_{\mathbf u_1^{k-1}}(a_{i,j})) - 1_{I_n} (T_{\mathbf u_1^{k-1}}(b_{i,j})) \big).\]
From Lemma~\ref{l:i4i5} and the symmetry of the map we get
\begin{equation}\nonumber
\KI_3(1)- \KI_3(1-\alpha)= 0 = \KI_3(-1)- \KI_3(\alpha-1),
\end{equation}
implying that $A_{3,1}=A_{3,4}=0$, $A_{3,2}=-1$ and $A_{3,3}=1$. Hence, any solution vector $\hat \gamma$ for $A \hat \gamma=\mathbf 0$ has the form $\hat \gamma= (\hat \gamma_1, \hat \gamma_2, \hat \gamma_2, \hat \gamma_3)$ and 
\eqref{e:fdensity} becomes
\begin{equation}\label{q:fdensity2}
 \begin{split}
f_p =  (\gamma_1 & + \gamma_2) \frac{p_1}{2}\sum_{k =0}^{M-2} \sum_{\mathbf u \in \Omega^k} \frac{p_{\mathbf u}}{2^k} \big( 1_{[-1, T_\mathbf u(\alpha-1))} - 1_{[-1,T_\mathbf u(-1))} \big)\\
& + (\gamma_2 + \gamma_3)\frac{p_0}{2} \sum_{k = 0}^{M-2} \sum_{\mathbf u \in \Omega^k} \frac{p_{\mathbf u}}{2^k} \big( 1_{[-1, T_\mathbf u(1))} - 1_{[-1,T_\mathbf u(1-\alpha))} \big),
\end{split}\end{equation}
where $\gamma = (\gamma_1, \gamma_2, \gamma_2, \gamma_3)$ is the unique non-trivial vector in the null space of the fundamental matrix $A$ that makes $f_p$ into a probability density function. In the next section we will derive a number of properties of $f_p$ with the goal of determining the frequency of the digit 0 in the signed binary expansions of $m_p \times \mu_p$ typical points.


\subsection{The frequency of the digit 0 in random signed binary expansions}
Recall from \eqref{q:randomdigits} that the random signed binary expansion of a point $(\omega,x)$ has a digit 0 in the $n$-th position precisely if
\[ R^{n-1}(\omega,x) \in [1]\times I_2 \cup \Omega^\mathbb N \times I_3 \cup [0] \times I_4 =: D_0.\]
Since $R$ is ergodic with respect to $m_p \times \mu_p$, it follows from Birkhoff's Ergodic Theorem that the frequency of the digit 0 in $m_p \times \mu_p$-almost all $(\omega,x)$ equals
\begin{equation}\label{q:freq0}
\pi_0 (\alpha, p) :=  \lim_{n \to \infty} \frac1n \sum_{k=0}^{n-1} 1_{D_0} (R^k(\omega,x)) = (1-p) \mu_p (I_2) +  \mu_p (I_3) + p \mu_p ( I_4 ).
\end{equation}
To give an example, consider $\alpha=1$, see Figure~\ref{f:rsb}(a). It is straightforward to check that the probability density $f_p = (1-p)1_{[-1,0]} + p1_{[0,1]}$ is invariant. This gives
\begin{equation}\label{q:1freq0}
\pi_0(1, p) = p \mu_p \Big( \Big[0, \frac12 \Big] \Big) + (1-p) \mu_p \Big( \Big[-\frac12,0 \Big] \Big) = \frac{p^2 + (1-p)^2}{2} \le \frac12
\end{equation}
with equality only for $p=0$ or $p=1$.

\vskip .2cm
To estimate $\pi_0 (\alpha, p)$ for other values of $\alpha$ we use a few lemmata. For $k \ge 1$ set $E_k = \{ \mathbf u \in \Omega^k \, : \, T_\mathbf u(1)=S^k(1)\}$ and $F_k =  \{ \mathbf u \in \Omega^k \, : \, T_\mathbf u(1-\alpha)=S^k(1)\}$. Also, use $(b_n)_{n \ge 1}$ to denote the digits in the signed binary expansion of 1 generated by $S$, i.e.,
\[ b_n = \begin{cases}
-1, & \text{if } S^{n-1}(1) < -\frac12,\\
0, & \text{if } -\frac12 \le S^{n-1}(1) \le \frac12,\\
1, & \text{if } S^{n-1}(1) > \frac12.
\end{cases}\]
Write $\mathbf b_k = b_1 \cdots b_k$ for any $k \ge 1$.

\begin{lem}\label{l:ekfk}
For all $1 \le k < M-1$, $F_k \subseteq E_k$ and $E_k \setminus F_k = \{ \mathbf b_k \}$.
\end{lem}

\begin{proof}
First note that the $n$-th signed binary digit of 1 generated by $S$, $n < M-1$, equals $0$ if $S^{n-1}(1) \in I_4$ and $1$ if $S^{n-1}(1) \in I_5$. We prove the statement by induction. For $k=1$ we have $E_1 = \{ 0,1\}$ and $F_1 = \{0\}$. Assume the statement holds for some $1 \le k <M-2$. If $S^k(1) = T_{\mathbf b_k}(1) \in I_4$, then $b_{k+1}=0$ and we know from the assumptions and since $S^k(1)-\alpha \in I_1$ that
\[ T_{\mathbf b_k0}(1) = S^{k+1}(1), \, T_{\mathbf b_k1}(1)= T_{\mathbf b_k0} (1-\alpha) = T_{\mathbf b_k1}(1-\alpha) = S^{k+1}(1)-\alpha.\]
Hence, $\mathbf b_k0 \in E_{k+1} \setminus F_{k+1}$ and $\mathbf b_k1 \not \in E_{k+1}\cup F_{k+1}$. If $S^k(1) = T_{\mathbf b_k}(1) \in I_5$, then $b_{k+1}=1$ and
\[ T_{\mathbf b_k0} (1) = T_{\mathbf b_k1}(1) = T_{\mathbf b_k0} (1-\alpha) = S^{k+1}(1), \,  T_{\mathbf b_k1}(1-\alpha) = S^{k+1}(1)-\alpha.\]
So, $\mathbf b_k1 \in E_{k+1} \setminus F_{k+1}$ and $\mathbf b_k0 \in E_{k+1} \cap F_{k+1}$. For any other $\mathbf u \in \Omega^k$ it holds that $T_\mathbf u(1) = T_\mathbf u(1-\alpha)$, so that either $\mathbf u j \in E_{k+1} \cap F_{k+1}$ or $\mathbf u j \not \in E_{k+1} \cup F_{k+1}$, $j=0,1$. This gives the statement.
\end{proof}

\begin{lem}\label{l:constant}
The density $f_p$ is constant and equal to $\frac1{\alpha}$ on the interval $[1-\alpha, \alpha-1]$.
\end{lem}

\begin{proof}
For any $\mathbf u \in \Omega^k$, write $\bar{\mathbf u} = (1-u_1)\cdots (1-u_k)$ and for a subset $E \subseteq \Omega^k$ write $\bar E = \{ \mathbf u \in \Omega^k \, : \, \bar{\mathbf u} \in E\}$.  
By Lemma~\ref{l:ekfk} we have for each $k < M$,
\begin{equation}\nonumber
\delta_k:= \sum_{ \mathbf u \in  E_k}  \frac{p_\mathbf u}{2^k}- \sum_{  \mathbf u \in F_k}  \frac{p_\mathbf u}{2^k} = \frac{p_{\mathbf b_k}}{2^k},\quad  \bar \delta_k := \sum_{\mathbf u \in  \bar E_k^c} \frac{p_\mathbf u}{2^k}- \sum_{ \mathbf u \in  \bar F_k^c} \frac{p_\mathbf u}{2^k} = \frac{p_{\bar{\mathbf b}_k}}{2^k},
\end{equation}
Recall the formula for the density $f_p$ from \eqref{q:fdensity2}. Using Proposition~\ref{p:distancealpha} we get
\[ \begin{split}
 \frac{p_0}{2}\sum_{k=0}^{M-2} & \sum_{\mathbf u \in \Omega^k}  \frac{p_\mathbf u}{2^k} \big( 1_{[-1, T_\mathbf u (1))} - 1_{[-1, T_\mathbf u(1-\alpha))} \big) \\
=\ &  \frac{p_0}{2}  \sum_{k=0}^{M-2} \sum_{\substack{\mathbf u \in \Omega^k:\\ T_\mathbf u(1)=S^k(1),\\ T_\mathbf u(1-\alpha)=S^k(1)-\alpha}} \frac{p_\mathbf u}{2^k} 1_{[T_\mathbf u(1-\alpha), T_\mathbf u (1))} - \frac{p_0}{2}\sum_{k=0}^{M-2} \sum_{\substack{\mathbf u \in \Omega^k:\\T_\mathbf u(1)=S^k(1)-\alpha, \\T_\mathbf u(1-\alpha)=S^k(1)}} \frac{p_\mathbf u}{2^k} 1_{[T_\mathbf u (1),T_\mathbf u(1-\alpha))}\\
=\ &  \frac{p_0}{2}  \sum_{k=0}^{M-2} \delta_k 1_{[S^k(1)-\alpha, S^k(1))}.
\end{split}\]
For the other side it holds similarly using the symmetry of the system that
\[  \frac{p_1}{2} \sum_{k=0}^{M-2} \sum_{\mathbf u \in \Omega^k} \frac{p_\mathbf u}{2^k} \big( 1_{[-1, T_\mathbf u (\alpha-1))} - 1_{[-1, T_\mathbf u(-1))} \big) =  \frac{p_1}{2}  \sum_{k=0}^{M-2}  \bar \delta_k 1_{[-S^k(1), \alpha- S^k(1))}. \]
By \eqref{q:distancealpha} we have for all $k < M-1$,
\[ S^k(1), \alpha - S^k(1) \in [\alpha-1, 1] \quad \text{ and } \quad S^k(-1), S^k(\alpha-1) \in [-1, 1-\alpha],\]
so that on $[1-\alpha, \alpha-1]$ we obtain
\begin{equation}\nonumber
f_p \mid_{[1-\alpha, \alpha-1]} (x)= (\gamma_1+\gamma_2) \frac{p_1}{2} \sum_{k=0}^{M-2}\bar \delta_k + (\gamma_2+\gamma_3) \frac{p_0}{2} \sum_{k=0}^{M-2} \delta_k.
\end{equation}
Since $f_p$ is a probability density it follows that
\begin{equation}\label{e:c}
1= \int_{[-1,1]} f_p \,  d\lambda = (\gamma_1+\gamma_2) \frac{p_1}{2}  \sum_{k=0}^{M-2} \bar \delta_k \alpha  + (\gamma_2+\gamma_3) \frac{p_0}{2} \sum_{k=0}^{M-2} \delta_k \alpha.
\end{equation}
Hence,
\begin{equation}\nonumber
f_p \mid_{[1-\alpha, \alpha-1]} (x) = \frac1\alpha,
\end{equation}
which gives the result.
\end{proof}

With this information we can compute $\pi_{0}(\alpha,p)$ for $\alpha \in \big[ \frac32,2\big]$. Since in this case $\alpha-1 \geq \frac12$ it follows from Lemma~\ref{l:constant} that
\begin{equation}\label{q:freq32}
\pi_{0}(\alpha, p) = \frac{\alpha-1}{\alpha} + \frac{2-\alpha}{2\alpha}= \frac12.
\end{equation}
That is, for $\alpha \geq \frac32$, and any $0 < p < 1$, the frequency of the digit 0 is equal to $\frac12$ in the signed binary expansion of $m_p \times \mu_p$-almost all $(\omega,x)$. For the other values of $\alpha$ we need to do some more work.

\begin{lem}\label{l:gamma2}
Let $\gamma = (\gamma_1, \gamma_2, \gamma_2, \gamma_3)$ be the unique vector in the null space of $A$ that makes $f_p$ into a probability density function. Then $\gamma_2 = \frac1{\alpha}$.
\end{lem}

\begin{proof}
Since $S^k(1) \in I_4 \cup I_5$ for all $k < M-1$ it follows from the definition of the function $\KI_n$ in \eqref{q:kin} and Proposition~\ref{p:distancealpha} that for $y=1, 1-\alpha$,
\[ \KI_4(y) + \KI_5(y) = \sum_{k=0}^{M-2} \sum_{j \in \Omega} \sum_{\mathbf u \in \Omega^k} \frac{p_j}{2} \frac{p_\mathbf u}{2^k} 1_{I_4 \cup I_5} (T_{\mathbf u}(y)) = \frac12 \sum_{k=0}^{M-2} \sum_{\stackrel{\mathbf u \in \Omega^k:}{T_{\mathbf u} (y)=S^k(1)}}  \frac{p_\mathbf u}{2^k},\]
so that
\[ \frac12 \sum_{k=0}^{M-2} \delta_k = \KI_4(1) - \KI_4 (1-\alpha) + \KI_5(1) -\KI_5 (1-\alpha).\]
A similar statement holds for $-1$ and $\alpha-1$. The fourth and fifth line of the linear system $A\gamma=\mathbf 0$ read
\[ p_1 (\KI_4(\alpha-1)-\KI_4(-1))(\gamma_1+\gamma_2) + p_0 (\KI_4(1)-\KI_4(1-\alpha))(\gamma_2+\gamma_3) - \gamma_2 +\gamma_3 =0\]
and
\[ p_1 (\KI_5(\alpha-1)-\KI_5(-1))(\gamma_1+\gamma_2) + p_0 (\KI_5(1)-\KI_5(1-\alpha))(\gamma_2+\gamma_3) - \gamma_3 =0,\]
respectively. Adding them up gives
\[ \begin{split}
\gamma_2 =\ & p_1(\gamma_1 + \gamma_2) (\KI_4 (\alpha-1)-\KI_4 (-1) + \KI_5 (\alpha-1)-\KI_5(-1))\\
& + p_0 (\gamma_2+\gamma_3) (\KI_4(1)-\KI_4 (1-\alpha) + \KI_5 (1)-\KI_5(1-\alpha))\\
=\ &  \frac{p_1}{2} (\gamma_1 + \gamma_2) \sum_{k=0}^{M-2} \bar \delta_k + \frac{p_0}{2} (\gamma_2 + \gamma_3) \sum_{k=0}^{M-2} \delta_k.
\end{split}\]
The result then follows from \eqref{e:c}.
\end{proof}

Combining Lemma~\ref{l:constant} and Lemma~\ref{l:gamma2} gives the following expression for the density $f_p$:
\begin{equation}\label{q:densityagain}
f_p =  \Big(\gamma_1 + \frac1{\alpha}\Big) \frac{p_1}{2}\sum_{k =0}^{M-2} \frac{p_{\bar{\mathbf b}_k}}{2^k} 1_{[-S^k(1), \alpha- S^k(1))}  + \Big(\frac1{\alpha} + \gamma_3 \Big)\frac{p_0}{2} \sum_{k = 0}^{M-2} \frac{p_{\mathbf b_k}}{2^k}1_{[S^k(1)-\alpha, S^k(1))},
\end{equation}
where $\mathbf b_k = b_1 \cdots b_k$ denote the first $k$ digits in the signed binary expansion of 1 given by $S$.

\begin{lem}\label{l:equivalent}
Let $\alpha \in \big(1, \frac32 \big)$ be a parameter for which the random system $R$ has strong random matching. Then both $\gamma_1, \gamma_3 \ge 0$. As a consequence, $f_p>0$ and $\mu_p$ is equivalent to the Lebesgue measure.
\end{lem}

\begin{proof}
Let $\gamma = (\gamma_1, \gamma_2, \gamma_2, \gamma_3)$ be the unique vector in the null space of $A$ that makes $f_p$ into a probability density. Set 
\[y = \max_{k \in \{1,2, \ldots, M-2\}} \{S^k(1), \alpha- S^k(1) \}.\]
By Lemma~\ref{l:difforbit} we can assume that $y \neq 1$. Then
\[ \mu_p ([y,1]) = m_p \times \mu_p (R^{-1}(\Omega \times [y,1])) = p \mu_p \Big( \Big[\frac{y-\alpha}{2}, \frac{1-\alpha}{2} \Big] \cup \Big[ \frac{y}{2}, \frac12 \Big] \Big).\]
By the definition of $y$ one can see from \eqref{q:densityagain} that
\[ \mu_p ([y,1]) = \frac{p(\gamma_2+\gamma_3)}{2}(1-y).\]
Furthermore, $T_0(1-\alpha)= 2-\alpha$ and $T_1(1-\alpha)= 2-2\alpha$, so in particular $y \geq \max \{2-\alpha, 2\alpha-2\}$. It follows that $1-\alpha \le \frac{y-\alpha}{2} < \frac{1-\alpha}{2}$. Thus by Lemma \ref{l:constant} and Lemma~\ref{l:gamma2}, $f_p\mid_{[\frac{y-\alpha}{2}, \frac{1-\alpha}{2}]}=\gamma_2$. We proceed by showing that none of the points $S^k(1)$ or $\alpha-S^k(1)$, $1 \le k \le M-2$, lie in the interval $\big[\frac{y}{2}, \frac12 \big]$, which then by \eqref{q:densityagain} implies that the density $f_p$ is also constant on the interval $\big[\frac{y}{2}, \frac12 \big]$. For $k= M-2= M_\alpha-1$, matching for $S$ implies that $\frac12 < S^{M-2}(1) <\alpha - \frac12$ and $ \alpha - S^{M-2}(1)  > \frac12$. Suppose there exists a $k \in \{1,2, \ldots, M-3\}$ such that $ \frac{y}{2} <  S^k(1) < \frac12$ (or $ \frac{y}{2} <  \alpha-S^k(1) < \frac12$). Then $ S^{k+1}(1) > y$ (or $ \alpha-S^{k+1}(1) > y$), which gives a contradiction with the definition of $y$. The same holds for $\alpha-S^k(1)$. Hence, there is a constant $c \ge 0$ such that 
\[ \frac{p(\gamma_2+\gamma_3)}{2}(1-y) =  \mu_p ([y,1]) = p \Big( \gamma_2 \frac{(1-y)}{2} + c \frac{(1-y)}{2} \Big).\]
So, $0 \le \mu_p \big( \big[ \frac{y}{2}, \frac12 \big] \big) = c=\gamma_3$. The proof that $\gamma_1\ge 0$ goes similarly. The fact that $f_p$ is strictly positive and the equivalence of $\mu_p$ and the Lebesgue measure now follow from \eqref{q:densityagain}.
\end{proof}

The following result can be proven in essentially the same way as \cite[Theorem 4.1]{DK}. We include a proof here for convenience.
\begin{lem}[cf.~Theorem 4.1 of \cite{DK}]\label{l:ct}
Fix $0 < p < 1$. The map $\alpha \mapsto \pi_0 (\alpha, p)$ is continuous on $\big( 1,\frac32 \big)$.
\end{lem}

\begin{proof}
In this proof we use $f_\alpha = f_{\alpha,p}$ to denote the unique density from \eqref{q:densityagain}. By \eqref{q:freq0}, for the continuity of $\alpha \mapsto \pi_0 (\alpha, p)$ it is sufficient to prove $L^1$-convergence of the densities $f_{\alpha}$; i.e., for any sequence $\{\alpha_k\}_{k \geq 1} \subseteq \big(1,\frac32 \big)$ converging to a fixed $\hat{\alpha} \in \big(1,\frac32\big)$, there is convergence $f_{\alpha_k} \rightarrow f_{\hat{\alpha}}$ in $L^1(\lambda)$. The proof of this fact goes along the following lines:
\begin{itemize}
\item[1.] First we show that there is a uniform bound, i.e., independent of $k$, on the total variation and supremum norm of the densities $f_{\alpha_k}$. It then follows from Helly's Selection Theorem that there is some subsequence of $(f_{\alpha_k})$ for which an a.e.~and $L^1$ limit $\hat f$ exist.
\item[2.] We show that $\hat f = f_{\hat \alpha}$, which by the same proof implies that any subsequence of $(f_{\alpha_k})$ has a further subsequence converging a.e.~to the same limit $f_{\hat \alpha}$. Hence, $(f_{\alpha_k})$ converges to $f_{\hat \alpha}$ in measure.
\item[3.] By the uniform integrability of $(f_{\alpha_k})$ it then follows from Vitali's Convergence Theorem that the convergence of $(f_{\alpha_k})$ to $f_{\hat \alpha}$ is in $L^1$.
\end{itemize}

Step 1.~and 2.~use Perron-Frobenius operators. For $j=0,1$ the Perron-Frobenius operator $P_{\alpha,j}$ of $T_{\alpha,j}$ is uniquely defined by the equation
\[ \int (P_{\alpha,j} f) g \, d\lambda = \int f (g \circ T_{\alpha,j}) \, d\lambda \quad \forall f \in L^1(\lambda), \, g \in L^\infty(\lambda)\]
and the Perron-Frobenius operator $P_\alpha$ of $R_\alpha$ is then defined by $P_{\alpha} f = p P_{\alpha,0} f + (1-p) P_{\alpha,1}f$.
Equivalently, $P_\alpha$ is uniquely defined by the equation
\begin{equation}\label{q:randompf}
\int (P_\alpha f)g \, d\lambda = p \int f (g \circ T_{\alpha,0} ) \, d\lambda + (1-p) \int f(g \circ T_{\alpha,1}) \, d\lambda \quad \forall f \in L^1(\lambda), \, g \in L^\infty(\lambda).
\end{equation}
Since each $R_\alpha$ has a unique probability density $f_\alpha$ it follows from \cite[Theorem 1]{Pe} that $f_\alpha$ is the $L^1$ limit of $(\frac1n \sum_{j=0}^{n-1} P_\alpha^j 1)_{n \geq 1}$ and it is the unique probability density that satisfies $P_\alpha f_\alpha = f_\alpha$.  From \cite[Theorem 5.2]{In} each $f_\alpha$ is a function of bounded variation. We proceed by finding uniform bounds on the total variation and supremum norm of these densities.

\vskip .2cm
Fix $\hat{\alpha} \in \big(1,\frac32\big)$. For the second iterates of the Perron-Frobenius operators we have
\[ P_{\alpha}^2 f = \sum_{i,j =0}^1 p_ip_j P_{\alpha,j} ( P_{\alpha,i} f).\]
Since the intervals of monotonicity of any of the maps $T_{\alpha,\mathbf u}$ for $\mathbf u \in \Omega^2$, only become arbitrarily small for $\alpha$ approaching $1$ and $\frac32$, we can find a uniform lower bound $\delta$ on the length of the intervals of monotonicity of any map $T_{\alpha,\mathbf u}$, $\mathbf u \in \Omega^2$, for all values $\alpha$ that are close enough to $\hat \alpha$. Applying \cite[Lemma 5.2.1]{BoGo} to $T_{\alpha,j}$, $j=0,1$, and any of the second iterates $T_{\alpha, \mathbf u}$, $\mathbf u \in \Omega^2$, gives that
\[ Var (P_{\alpha,j} f) \le Var (f) + \frac1{\delta} \|f\|_1 \quad \text{ and } \quad Var (P_{\alpha,\mathbf u}f) \le \frac12 Var(f) + \frac1{2\delta}\|f \|_1,\]
where $Var$ denotes the total variation over the interval $[-1,1]$. Since these bounds do not depend on $\alpha, j , \mathbf u$, the same estimates hold for $P_\alpha$, so that for any function $f: [-1,1] \to \mathbb R$ of bounded variation and any $n \ge 1$,
\begin{equation}\label{q:lyn}
Var (P_\alpha^n f) \le \frac1{2^{\lfloor n/2 \rfloor}}Var (f)  + \Big(2+\frac1{\delta} \Big) \|f \|_1.
\end{equation}
Let $\{\alpha_k\}_{k \geq 1}$ with $\alpha_k \rightarrow \hat{\alpha}$ be a sequence for which the lower bound $\delta$ holds for each $k$. For each $k$ and $n$, write $f_{k,n} = \frac1n \sum_{i=0}^{n-1} P_{\alpha_k} 1$. Since
\[ \sup |f_{k,n}| \le Var (f_{k,n}) + \int f_{k,n} \, d\lambda,\]
it follows from \eqref{q:lyn} that there is a uniform constant $C>0$ (independent of $k,n$) such that $Var (f_{k,n}), \sup |f_{k,n}| < C$. The same then holds for the limits $f_{\alpha_k}$. Helly's Selection Theorem then gives the existence of a subsequence $\{k_i\}$ and a function $\hat f$ of bounded variation, such that $f_{\alpha_{k_i}} \to \hat f$ in $L^1(\lambda)$ and $\lambda$-a.e.~and with $Var(\hat f), \sup |\hat f| < C$. This finishes 1.

\vskip .2cm
By 2.~and 3.~above, what remains to finish the proof is to show that $P_{\hat \alpha} \hat f = \hat f$. By \eqref{q:randompf} it is enough to show that for any compactly supported $C^1$ function $g:[-1,1]\to \mathbb R$ it holds that
\[ \left| \int (P_{\hat \alpha} \hat f) g \, d\lambda - \int \hat f g \, d\lambda \right|=0.\]
Note that
\[ \begin{split} 
\left| \int (P_{\hat \alpha} \hat f) g \, d\lambda - \int \hat f g \, d\lambda \right| \le p \left| \int \hat f (g \circ T_{\hat \alpha,0}) \, d\lambda - \int \hat f g \, d\lambda \right| + (1-p) \left| \int \hat f (g \circ T_{\hat \alpha,1}) \, d\lambda - \int \hat f g \, d\lambda \right|.
\end{split} \]
For $j=0,1$ we can write
\[\begin{split}
\left| \int \hat f (g \circ T_{\hat \alpha,j}) \, d\lambda - \int \hat f g \, d\lambda \right| \le \ & \left| \int \hat f (g \circ T_{\hat \alpha,j}) \, d\lambda - \int f_{\alpha_{k_i}} (g \circ T_{\hat \alpha,j}) \, d\lambda \right| \\
& + \left| \int f_{\alpha_{k_i}} (g \circ T_{\hat \alpha,j}) \, d\lambda - \int f_{\alpha_{k_i}} (g \circ T_{ \alpha_{k_i},j}) \, d\lambda \right| \\
& + \left| \int f_{\alpha_{k_i}} (g \circ T_{\alpha_{k_i},j}) \, d\lambda - \int \hat f  g \, d\lambda \right| .
\end{split}\]
The first and third integral on the right hand side can be bounded by $\| g \|_\infty \| \hat f - f_{\alpha_{k_i}}\|_1 \to 0$. For the second integral, $\| f_{\alpha_{k_i}}\|_\infty < C$ and $\int |g \circ T_{\hat \alpha,j} - g \circ T_{ \alpha_{k_i},j}| \, d\lambda \to 0$ by the Dominated Convergence Theorem. Hence, $\hat f = f_{\hat \alpha}$ and $f_{\alpha_k}\to f_{\hat \alpha}$ in $L^1$.
\end{proof}

\begin{figure}[h]
\centering
\includegraphics[scale=0.3]{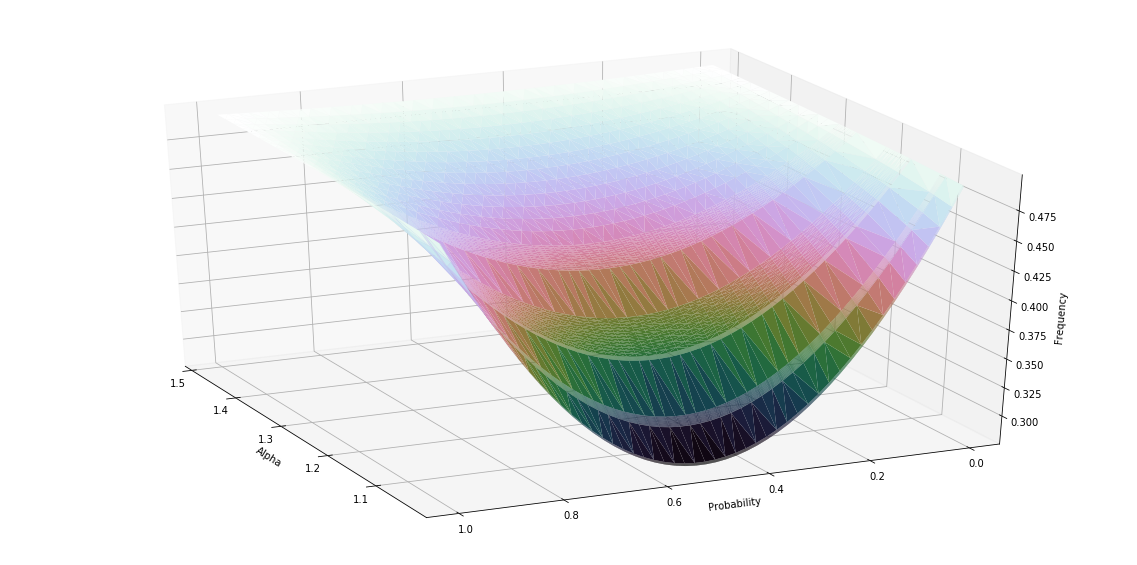}
\caption{The graph of $(\alpha,p) \mapsto \pi_{0}(\alpha,p)$.}
\label{f:nesimulation}
\end{figure}

Figure~\ref{f:nesimulation} shows a numerical approximation of the graph of the function $(\alpha,p) \mapsto \pi_{0}(\alpha,p)$. We can now prove that the maximal value of the frequency of the digit 0 is in fact $\frac12$.

\begin{thm}
For any $0 < p<1$ and any $\alpha \in [1,2]$ the frequency $\pi_0 (\alpha, p)$ is at most $\frac12$ for $m_p \times \lambda$-a.e.~$(\omega,x) \in \Omega^\mathbb N \times [-1,1]$.
\end{thm}

\begin{proof}
For $\alpha \in \big[ \frac32, 2 \big]$ the statement follows from \eqref{q:freq32} and for $\alpha=1$ from \eqref{q:1freq0}. Let $\alpha \in \big(1, \frac32 \big)$. The deterministic map $T_{\alpha,0}$ has density $f_0 = \frac1{\alpha}1_{[1-\alpha, 1]}$ and $T_{\alpha,1}$ has $f_1=\frac1{\alpha} 1_{[-1, \alpha-1]}$. Hence $\pi_0(\alpha, p) = \frac12$ for $p=0,1$. Let $0<p<1$ and let $\alpha$ be a parameter satisfying the conditions of Lemma~\ref{l:equivalent}. We know that $f_p$ is constant and equal to $\frac1{\alpha}$ on $[1-\alpha, \alpha-1]$. For $x > \alpha -1$ the density can be written as
\[ \begin{split}
f_p(x) =\ & \frac{1}{\alpha }- \bigg( \frac{(1-p)(\gamma_1 + \gamma_2)}{2} \sum_{k=0}^{M-2} \frac{p_{\bar{\mathbf b}_k}}{2^k} 1_{[\alpha-1, x ]} (\alpha-S^k(1))\\
&   + \frac{p(\gamma_2 + \gamma_3)}{2} \sum_{k=0}^{M-2} \frac{p_{\mathbf b_k}}{2^k} 1_{[\alpha-1,x] }(S^k(1)) \bigg)\\
=\ & \frac{1}{\alpha } - \frac{(1-p)(\gamma_1 + \gamma_2)}{2} - \bigg( \frac{(1-p)(\gamma_1 + \gamma_2)}{2} \sum_{k=1}^{M-2} \frac{p_{\bar{\mathbf b}_k}}{2^k} 1_{[\alpha-1, x ]} (\alpha-S^k(1)) \\
&+ \frac{p(\gamma_2 + \gamma_3)}{2} \sum_{k=1}^{M-2} \frac{p_{\mathbf b_k}}{2^k} 1_{[\alpha-1,x] }(S^k(1)) \bigg)\\
\le & \frac{1}{\alpha } - \frac{(1-p)(\gamma_1 + \gamma_2)}{2} .
\end{split}\]
Similarly, for $x < 1-\alpha$ we get $f_p(x) \le \frac{1}{\alpha } - \frac{p(\gamma_2 + \gamma_3)}{2}$. By \eqref{q:freq0} and Lemma~\ref{l:equivalent},
\begin{equation}\nonumber
\begin{split}
\pi_{0}(\alpha,p) =& \ \frac{\alpha-1}{\alpha} +  \frac{\alpha-1}{2\alpha} + p \mu_{\alpha, p}\Big(\Big[\alpha-1, \frac12\Big]\Big) + (1-p) \mu_{\alpha, p}\Big(\Big[-\frac12, 1-\alpha\Big]\Big)\\
\leq & \ \frac{3(\alpha-1)}{2\alpha} + \frac{3-2\alpha}{2\alpha} \bigg(1- \frac{p(1-p)\alpha}{2} \min\{\gamma_1+\gamma_2, \gamma_2+\gamma_3\} \bigg)\\
= & \ \frac12 - \frac{3-2\alpha}{2} \frac{p(1-p)}{2}  \min\{\gamma_1+\gamma_2, \gamma_2+\gamma_3\} \\
< & \ \frac12.
\end{split}
\end{equation}
Since matching holds for Lebesgue almost all parameters $\alpha$, the statement now follows from Lemma~\ref{l:ct} and the equivalence of $\mu_p$ and the Lebesgue measure.
\end{proof}

\section{Final remarks}

\subsection{Remarks on the symmetric doubling maps}
The numerical approximation of the graph of $(\alpha,p) \mapsto \pi_0(\alpha,p)$ shown in Figure~\ref{f:nesimulation} seems to suggest some other features of the map that we have not proved. Firstly, it suggests some symmetry. In fact it can be shown that for each fixed $\alpha$ and any $x \in [0,1]$, it holds that $f_p(x) = f_{1-p}(-x)$. For this one needs to consider the fundamental matrix $\tilde A$ corresponding to the random system $\tilde R_\alpha$ obtained by switching the roles of $p$ and $1-p$. Then using the permutation $(12)(45)$, one can relate various of the quantities involved for $\tilde A$ to the fundamental matrix $A$ of $R_\alpha$.

\vskip .2cm
Secondly, for any matching parameter $\alpha$ and any $0 < p < 1$ the density $f_{\alpha,p}$ is a finite combination of indicator functions, whose supports depend on the position of the points in the set $\{S^k(1), \alpha-S^k(1)\}_{k=0}^{M-2}$ and whose coefficients are polynomials in $p$. So, for such a fixed $\alpha$ and any $x \in [-1,1]$, the map $p \mapsto f_{\alpha,p}(x)$ is continuous in $p$. 

\vskip .2cm
Thirdly, the graph also suggests that the map presents a minimum at $p=\frac12$. Using the above two facts we were only able to show the following:
\begin{prop}
Let $\alpha \in [1,2]$ be such that $R$ has strong random matching. Then the map $p \mapsto \pi_0(\alpha,p)$ has an extremal value at $p=\frac12$. 
\end{prop}

\begin{proof}
By combining \eqref{q:freq0} and the fact that $f_p(x) = f_{1-p}(-x)$ we obtain
\[ \pi_0 (\alpha, p) = (1-p) \mu_{1-p} (I_4) +  \mu_p (I_3) + p \mu_p ( I_4 ).\]
Computing the derivative with respect to $p$ then gives
\begin{equation}\label{q:partial12}
\partial_p \pi_0(\alpha,p) = - \mu_{1-p} (I_4) -(1-p)\partial_p (\mu_{1-p} (I_4)) +  \partial_p(\mu_p (I_3)) +  \mu_p ( I_4 ) + \partial_p(\mu_p(I_4)).
\end{equation}
From Lemma~\ref{l:constant} it follows that $\partial_p(\mu_p (I_3)) = - \partial_p (\mu_{1-p}(I_3))$, implying that $\partial_p (\mu_p(I_3))=0$ at $p=\frac12$. Therefore, by \eqref{q:partial12} $\partial_p \pi_0(\alpha,p)=0$ at $p=\frac12$.
\end{proof}

\subsection{Remarks on random continued fractions}
Theorem~\ref{t:pwconstant} states that for random piecewise affine maps of the interval satisfying (c1), (c2) and (c3) strong random matching implies that there exists a piecewise constant invariant density. Condition \eqref{q:randomderivatives} was sufficient for the theorem to work, which was one of the main motivations for Definition~\ref{d:rmatching2}.

\vskip .2cm
Theorem~\ref{t:pwconstant} is a random analogue of \cite[Theorem 1.2]{BCMP}, except that there the statement has less assumptions. The authors mention in \cite[Remark 1.3]{BCMP} that for piecewise smooth interval maps with strong matching the corresponding invariant probability densities are piecewise smooth. On the other hand, as we noted before, the natural extension construction which for continued fraction transformations is often used to find invariant densities, seems to suggest that matching alone is sufficient to guarantee the existence of a piecewise smooth density. It would be interesting to investigate this further for the random continued fraction transformation.

\begin{figure}[h]
\centering
\subfigure[$\alpha = 0.70315\ldots$, $p_0=0.3$]{\includegraphics[width=.3\textwidth]{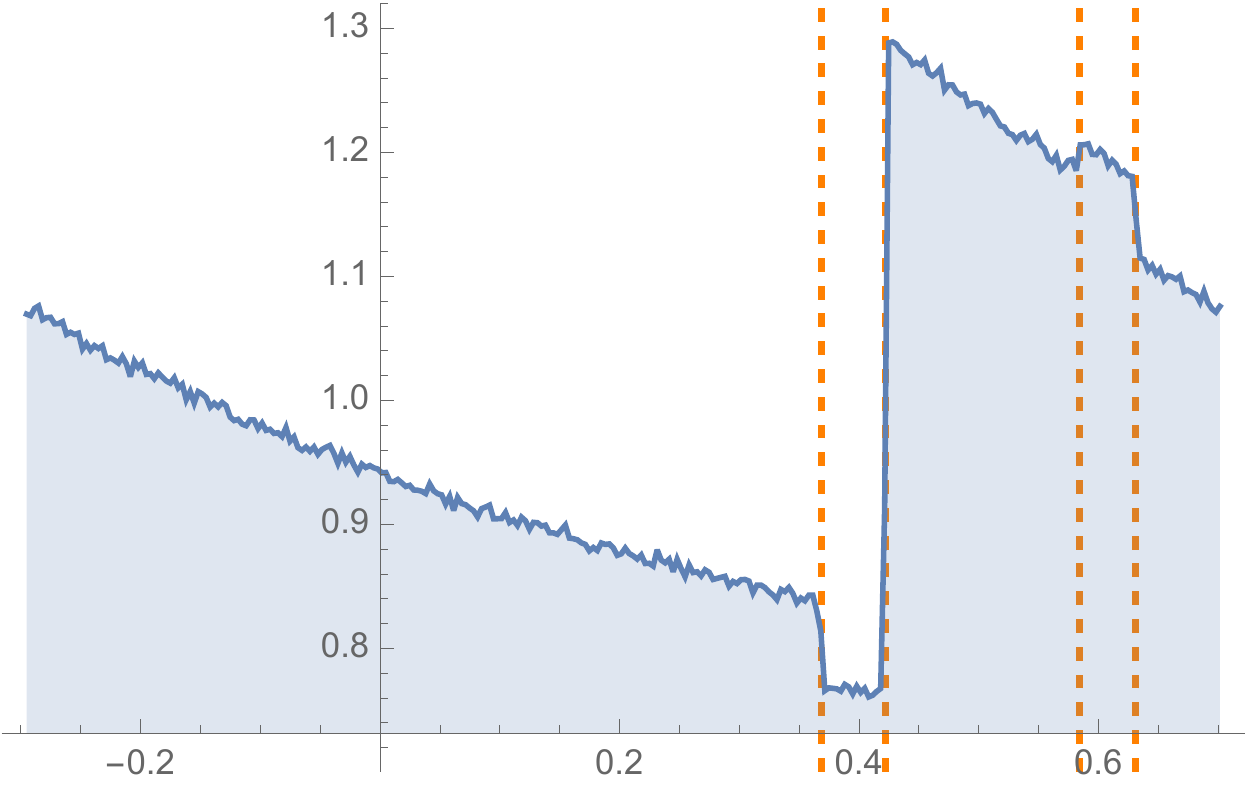}}
\hspace{1.5cm}
\subfigure[$\alpha = 0.77287\ldots$, $p_0=0.6$]{\includegraphics[width=.3\textwidth]{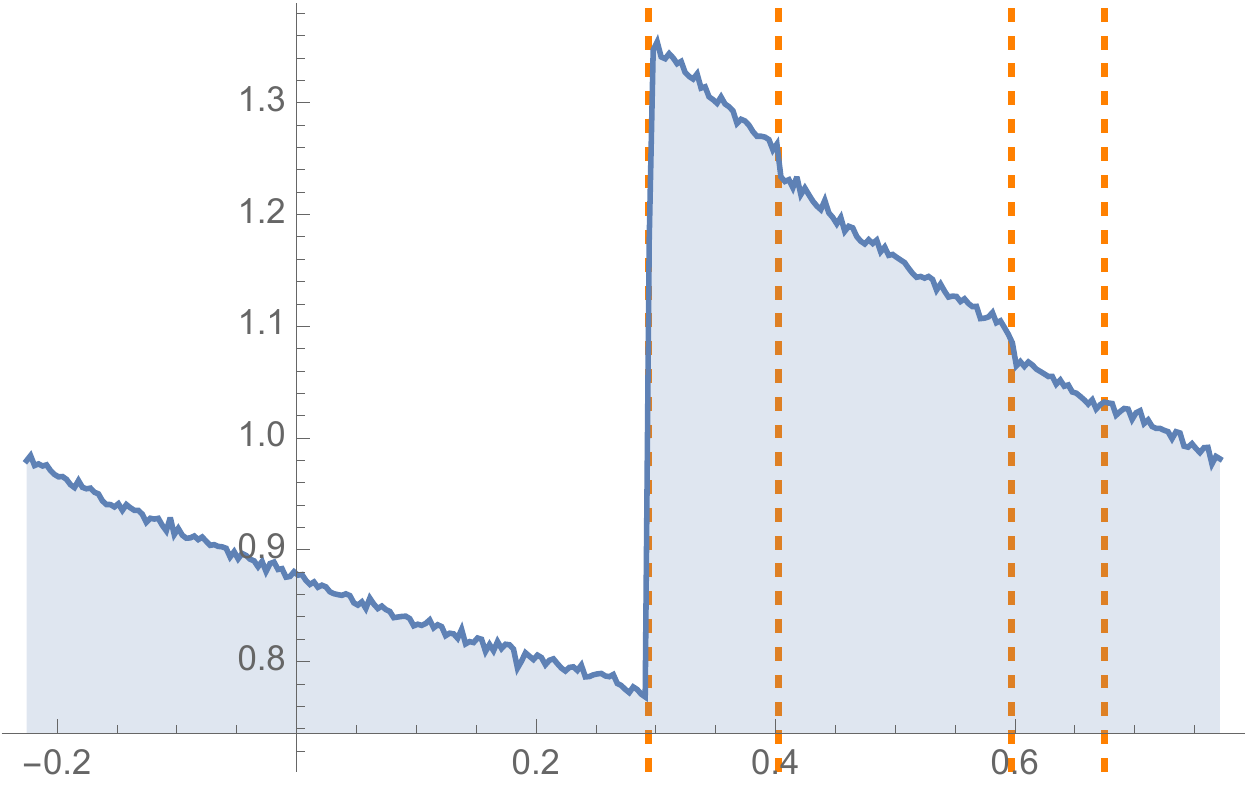}}
\caption{Numerical simulations of the invariant probability densities of the random continued fraction maps $R_\alpha$ from Example~\ref{x:randomcf}. In (a) we take $\alpha \in J_4$ and $p_0=0.3$ and in (b) we have $\alpha \in J_5$ and $p_0=0.6$. The dashed lines indicate the positions of the {\em prematching points}, i.e., the points in the orbits of $\alpha$ and $\alpha-1$ before the moment of matching.}
\label{f:sims}
\end{figure}

\vskip .2cm
In a first attempt to investigate to what extent Theorem~\ref{t:pwconstant} can be generalised to piecewise smooth random systems on an interval, we include some numerical simulations. Recall from Example~\ref{x:randomcf} that the random continued fraction maps $R_\alpha$ have strong random matching for $\alpha$ in the intervals $J_n$ with endpoints as in \eqref{q:jn}, see also Figure~\ref{f:intervalsIn}. Figure~\ref{f:sims} shows two simulations of the invariant densities for such systems $R_\alpha$. The densities seem to be piecewise smooth with discontinuities precisely at the orbit points of $\alpha$ and $\alpha-1$ before matching. This seems to support the claim that strong random matching is sufficient to guarantee the existence of a piecewise smooth invariant density.

\begin{figure}[h]
\centering
\subfigure[$\alpha = 0.584\ldots$, $p_0=0.3$]{\includegraphics[width=.25\textwidth]{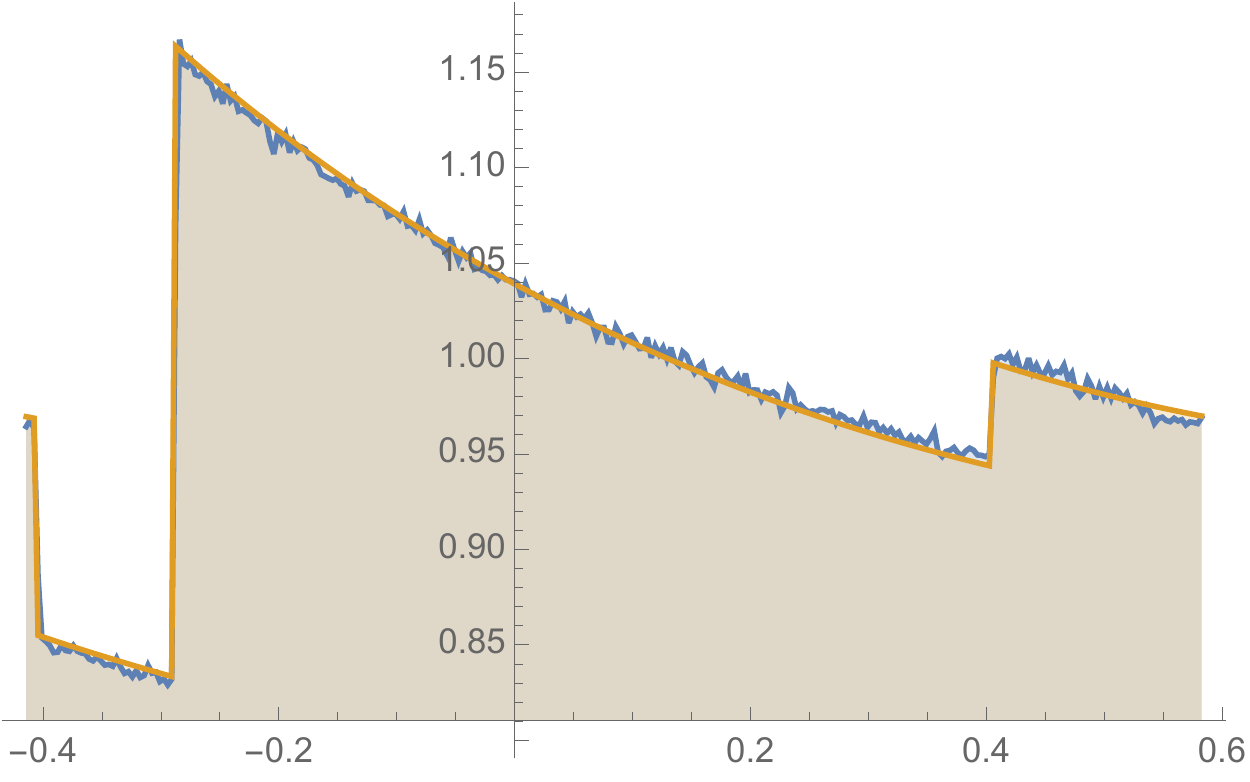}}
\hspace{1cm}
\subfigure[$\alpha = 0.579\ldots$, $p_0=0.65$]{\includegraphics[width=.25\textwidth]{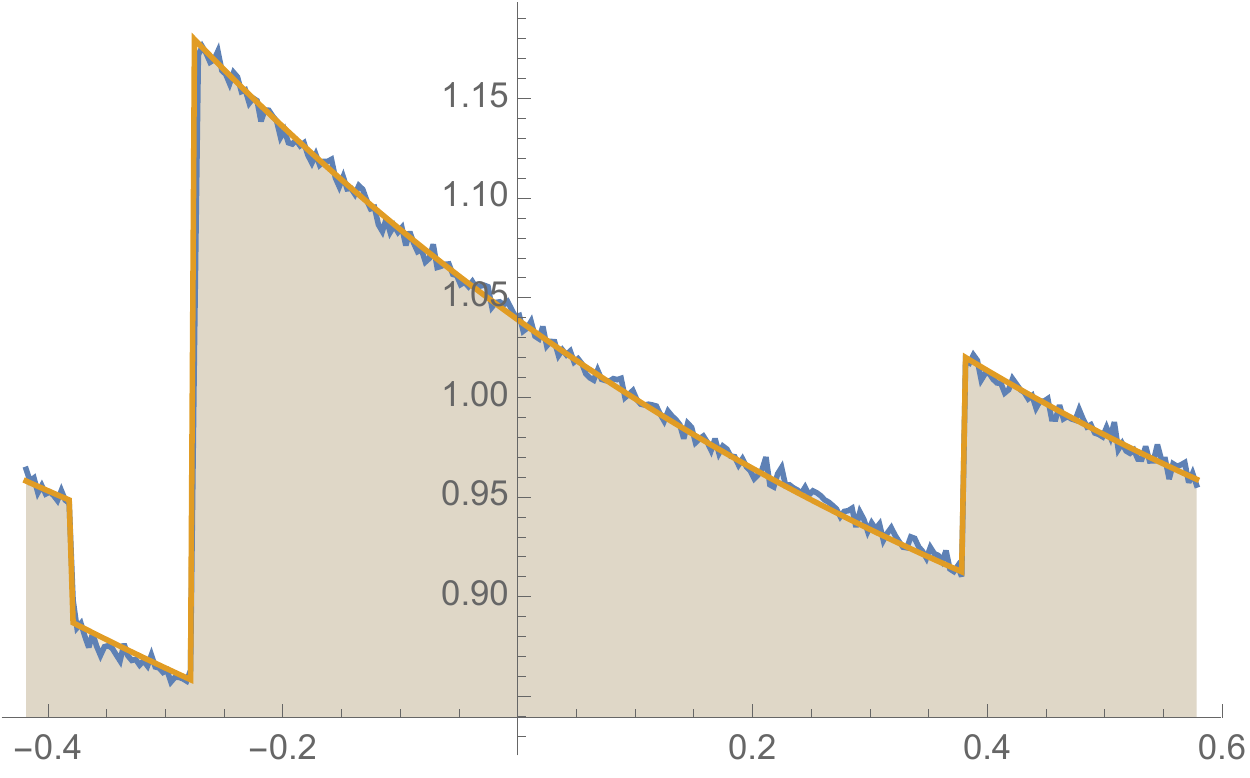}}
\hspace{1cm}
\subfigure[$\alpha = 0.541\ldots$, $p_0=0.25$]{\includegraphics[width=.25\textwidth]{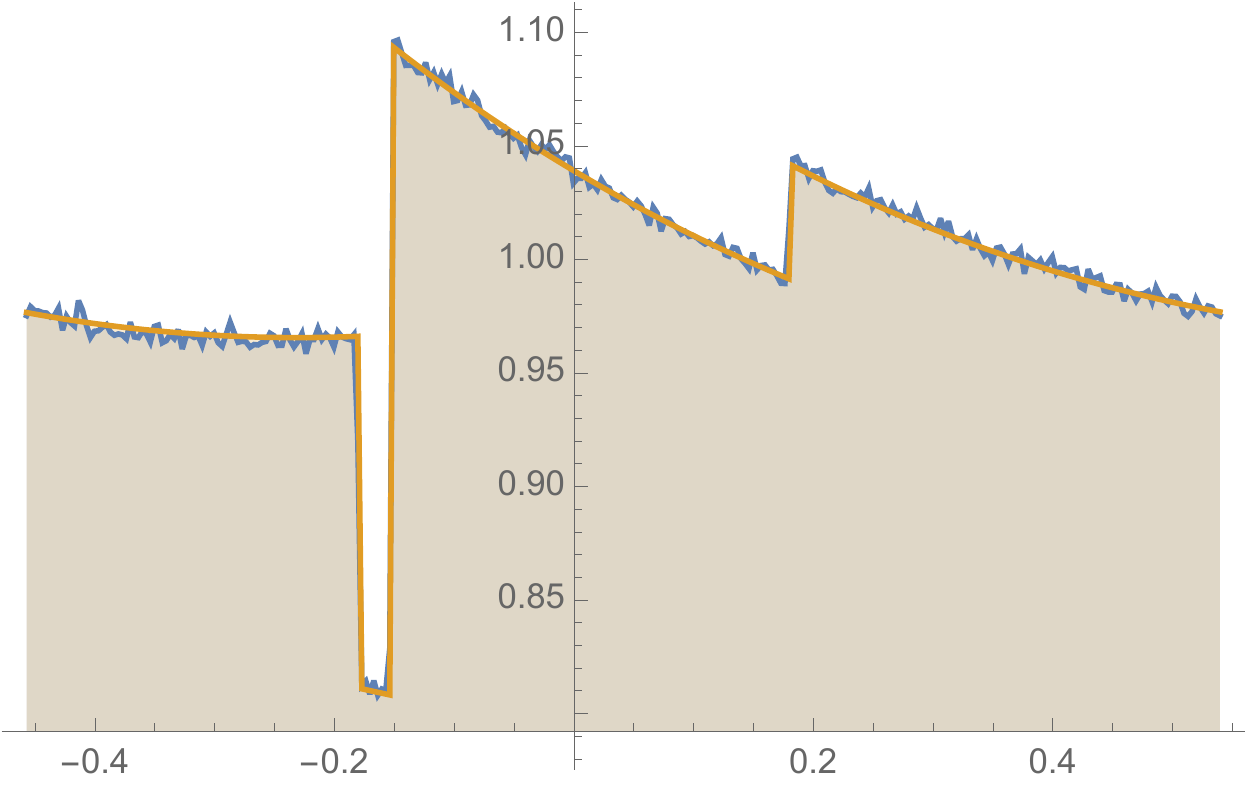}}
\caption{Numerical simulations of the invariant probability densities of the random continued fraction maps $R_\alpha$ from Example~\ref{x:randomcf} for three values of $\alpha$ between $\frac12$ and $2-\sqrt 2$. The map in $(a)$ has $\alpha \in \big( \frac{\sqrt{10}-2}{2}, 2-\sqrt 2 \big)$, which is the matching interval considered in Example~\ref{x:randomcf}. The orange graph is the graph of the weighted average of the densities of $T_{\alpha,0}$ and $T_{\alpha,1}$ with the appropriate values of $p$.}
\label{f:sims2}
\end{figure}

\vskip .2cm
In Example~\ref{x:randomcf} we also considered the maps $R_\alpha$ for $\alpha \in \big(\frac{\sqrt{10}-2}{2}, 2-\sqrt 2 \big)$. We showed that $R_\alpha$ has random matching with $M=3$, but no strong matching at that moment. With a similar approach it can be shown that $R_\alpha$ has random matching for various other intervals in $\big[ \frac12, \frac{\sqrt 5-1}{2} \big]$. For $\alpha \in \big[\frac12, 2-\sqrt 2\big]$ both deterministic maps $T_{\alpha,0}$ and $T_{\alpha,1}$ have strong matching with $M,Q \le 2$, as was shown in \cite{Nak81} and \cite{IT}, and moreover, for both of them the invariant densities are known. In Figure~\ref{f:sims2} we have plotted the weighted average of these densities together with numerical simulations of the densities for various values of $\alpha \in \big[ \frac12, 2-\sqrt{2} \big]$ and $0 < p < 1$. This makes us wonder whether we need strong random matching to guarantee the existence of a piecewise smooth invariant density for these random systems or whether random matching is sufficient.

\section{Acknowledgment}
The third author is supported by the NWO TOP-Grant No.~614.001.509.

\bibliographystyle{alpha}
\bibliography{random}

\end{document}